\documentclass[11pt]{article}

\usepackage[T1]{fontenc}
\usepackage{lmodern}
\usepackage[margin=1in]{geometry}
\usepackage{amsmath,amssymb,amsthm,mathtools,mathrsfs,comment}
\usepackage{microtype}
\usepackage[hidelinks]{hyperref}
\usepackage{enumitem}

\allowdisplaybreaks[2]
\newtheorem{theorem}{Theorem}[section]
\newtheorem{proposition}[theorem]{Proposition}
\newtheorem{lemma}[theorem]{Lemma}
\newtheorem{corollary}[theorem]{Corollary}
\theoremstyle{remark}
\newtheorem{remark}[theorem]{Remark}
\numberwithin{equation}{section}

\newcommand{\R}{\mathbb R}
\newcommand{\C}{\mathbb C}

\newcommand{\N}{\mathbb N}
\newcommand{\Q}{\mathbb Q}
\newcommand{\Z}{\mathbb Z}
\newcommand{\PP}{\mathbb P}
\newcommand{\EE}{\mathbb E}
\newcommand{\Md}{\mathfrak m_d}
\newcommand{\dd}{\,\mathrm d}
\newcommand{\ind}{\mathbf 1}

\newcommand{\Cont}{\operatorname{Cont}}
\newcommand{\dist}{\operatorname{dist}}
\newcommand{\diam}{\operatorname{diam}}

\newcommand{\dimH}{\dim_{\mathrm H}}
\newcommand{\dimP}{\dim_{\mathrm P}}
\newcommand{\dimBl}{\underline{\dim}_{\mathrm B}}
\newcommand{\dimBu}{\overline{\dim}_{\mathrm B}}
\newcommand{\Nc}{\mathcal N}

\newcommand{\HH}{\mathbb{H}}
\newcommand{\D}{\mathbb{D}}

\newcommand{\pa}{\partial}

\def\eps{\varepsilon}
\def\til{\widetilde}
\def\ha{\widehat}
\def\sem{\setminus}
\def\lin{\overline}
\def\ulin{\underline}
\def\vphi{\varphi}

\def\L{{\mathcal L}}

\newcommand{\BGE}{\begin{equation}}
\newcommand{\BGEN}{\begin{equation*}}
\newcommand{\EDE}{\end{equation}}
\newcommand{\EDEN}{\end{equation*}}

\title{Uniform Dimension Theorem for SLE in Natural Time}
\author{Dapeng Zhan}
\date{\today}

\begin{document}
\maketitle

\begin{abstract}
For $0<\kappa<8$, we prove a uniform dimension transformation
theorem for SLE in natural time: on a single probability-one event,
Hausdorff and packing dimensions of every time set are multiplied by
$d=1+\kappa/8$, and the same holds for lower and upper box dimensions
of every bounded time set.  The result applies to self-similar SLE
with stationary increments, chordal SLE in the half-plane, and
chordal SLE in bounded piecewise $C^1$ Jordan domains. In the last case, the box-dimension
conclusion holds under a positive-angle endpoint hypothesis.  We also
establish boundary-inclusive local H\"older continuity of every order
below $1/d$ for chordal SLE in natural time and identify the exact
local dimension $d$ of the natural occupation measure on compact time
intervals.
\end{abstract}



\section{Introduction}
Schramm--Loewner evolution (SLE), introduced by Schramm, is the
canonical family of conformally invariant random planar curves and
arises as the scaling limit of interfaces in a number of
two-dimensional critical lattice models; see
\cite{SchrammSLE,RohdeSchramm,LawlerBook} for its introduction,
basic properties, and further background.

Throughout this paper, we fix $\kappa\in(0,8)$ and write $d=1+\kappa/8\in(1,2)$,
which is the Hausdorff dimension of an SLE$_\kappa$ trace
\cite{Beffara}.  Lawler and Rezaei proved that the
$d$-dimensional Minkowski content of chordal SLE$_\kappa$ in $\HH$ from $0$ to $\infty$
exists and gives the natural parametrization \cite{LawlerRezaei}; the meaning of this
parametrization is recalled in \eqref{eq:intervalcontent}.
We use $\gamma$ to denote chordal SLE$_\kappa$ curves
with this parametrization.
A two-sided whole-plane SLE$_\kappa$ through the origin,
also parametrized by $d$-dimensional Minkowski content,
can be chosen as a self-similar process of index $1/d$
with stationary increments \cite{ZhanHolder};
see \eqref{ss} and \eqref{si} for the precise identities.
We denote this process by $\eta:\R\to\C$, normalized by
$\eta(0)=0$, and call it an \emph{sssi SLE$_\kappa$}.

For the sssi process $\eta$, it was proved in \cite[Theorem~1.3]{ZhanHolder}  that, for every deterministic closed set $A\subset\R$, $\dimH\eta(A)=d\dimH A$
almost surely.  We strengthen this to a uniform statement on a single
probability-one event, valid simultaneously for every subset of time,
including sets selected after the sample path is known.  We also obtain
the corresponding identities for packing dimension and, on bounded
time sets, for both box dimensions. We further establish
corresponding uniform dimension transformation formulas for chordal
SLE$_\kappa$ in the upper half-plane from $0$ to $\infty$, and for chordal
SLE$_\kappa$ in bounded piecewise $C^1$ Jordan domains.

The classical analogue is Kaufman's uniform dimension-doubling
theorem for Brownian motion; see \cite[Section~9.4]{MortersPeres}.
The inverse-image covering argument underlying Kaufman's theorem
extends naturally to packing dimension and, for bounded time sets,
to lower and upper box dimensions.  We formulate below a deterministic
version of this principle that treats all four notions of dimension
simultaneously.

We write $\HH=\{z\in\C:\Im z>0\}$ for the upper half-plane,
$\dimH$ and $\dimP$ for Hausdorff and packing dimensions,
and $\dimBl$ and $\dimBu$ for lower and upper box dimensions,
respectively. The dimension definitions will be recalled in
Section~\ref{sec:deterministic}.
Unless otherwise indicated, constants implicit in $\lesssim$ may
depend on $\kappa$, and subscripts indicate additional dependencies.
In some sections, further parameters are fixed and allowed to enter
the implicit constants; this will be stated explicitly there.
All time sets $A$ in the theorems below are arbitrary subsets
of the indicated time domains, with no measurability assumption.

\begin{theorem}\label{thm:main}
Let $\eta:\R\to\C$ be the sssi SLE$_\kappa$ process.
There exists an event of probability one on which, simultaneously
for every $A\subset\R$,
\begin{equation}\label{eq:main}
 \mathfrak D(\eta(A))=d\,\mathfrak D(A),
 \qquad
 \mathfrak D\in\{\dimH,\dimP,\dimBl,\dimBu\},
\end{equation}
where $A$ is required to be bounded for the two box dimensions.
\end{theorem}

\begin{theorem}\label{thm:chordal}
Let $\gamma:[0,\infty)\to\overline\HH$ be chordal SLE$_\kappa$
from $0$ to $\infty$ in $\HH$, parametrized by its
$d$-dimensional Minkowski content.
There exists an event of probability one on which, simultaneously
for every $A\subset[0,\infty)$, \eqref{eq:main} holds with $\gamma$ in place of $\eta$,
where $A$ is required to be bounded for the two box dimensions.
\end{theorem}

A \emph{regular $C^1$ arc} is the image of an injective $C^1$ map
$\varphi:[0,1]\to\C$ such that $\varphi'(t)\ne0$ for every
$t\in[0,1]$. Suppose that $\partial D$ is the union of finitely many regular
$C^1$ arcs.  We say that the \emph{endpoint condition} holds at
$x\in\partial D$ if either $x$  is not a junction point of these arcs,
or $x$ is a junction point and the interior angle of $D$ at $x$ is
positive.  The angle $2\pi$ is allowed.

\begin{theorem}\label{thm:smooth}
Let $D\subset\C$ be a bounded Jordan domain whose boundary is the union
of finitely many regular $C^1$ arcs, and let $a,b\in\partial D$ be
distinct.  Consider chordal SLE$_\kappa$ in $D$ from $a$ to $b$, viewed
modulo reparametrization.  Almost surely, it admits a parametrization
$\gamma_D:[0,L]\to\overline D$ by $d$-dimensional Minkowski content, in
the sense of \eqref{eq:intervalcontent}, with $0<L<\infty$,
$\gamma_D(0)=a$, and $\gamma_D(L)=b$.

There exists an event of probability one on which, simultaneously for
every $A\subset[0,L]$, \eqref{eq:main} holds with $\gamma_D$ in place of $\eta$.
For the two box dimensions, the endpoint condition is required to hold
at both $a$ and $b$.
\end{theorem}

When $0<\kappa\le4$, substantially stronger versions of
Theorem~\ref{thm:smooth} hold, with weaker boundary assumptions;
see Remark~\ref{rem:simple-bounded-domain}.

The dimension transformation results share a common deterministic
criterion.  The upper dimension bounds follow directly from
H\"older continuity.  For the lower bounds, the same H\"older control
is combined with uniform occupation estimates in an inverse-image
covering argument. Indeed, fix a compact time interval and exponents
$0<\alpha<1/d$ and $0<\beta <d$.
H\"older continuity ensures that each time-mesh interval of length
comparable to $r^{1/\alpha}$ which visits a ball of radius $r$
is mapped into the concentric ball of radius $2r$.
An occupation bound of order $r^\beta$ therefore limits the number
of such intervals to order $r^{\beta-1/\alpha}$.
Pulling back spatial covers yields the lower dimension bounds
as $\alpha\uparrow1/d$ and $\beta \uparrow d$.
Since these covering estimates hold before the time set is chosen,
the resulting identities hold simultaneously for every time set.

The H\"older input is different for the two processes.  For the sssi
SLE, the required regularity is already known from
\cite[Theorem~1.2]{ZhanHolder}.  For chordal SLE in $\HH$ from $0$ to
$\infty$, however, we need natural-time H\"older continuity up to the
boundary, including at the starting point and at later boundary
contacts.  Interior local equivalence alone does not provide this
control in the non-simple range.  We obtain the required
boundary-inclusive estimate by combining whole-plane estimates with a
boundary crossing estimate; see Proposition~\ref{prop:chordal-holder}.

The occupation estimate does not appear to be available in the
literature in the form needed here.  We first prove the required
estimate for chordal SLE in $\HH$, using the multipoint approach
estimate in \cite{RezaeiZhan} together with the Minkowski-content
results in \cite{LawlerRezaei}. The corresponding estimate for the sssi process is obtained by first
proving a multipoint approach estimate for its whole-plane
SLE$_\kappa(2)$ arm, using the conformal Markov property and a
conditioning argument for chordal SLE.  Combining these occupation bounds with the
corresponding H\"older estimates and the deterministic criterion
proves the first two dimension theorems.

For bounded piecewise $C^1$ Jordan domains, an additional issue is
that the whole curve, including its approach to the two marked
boundary points, must first be shown to have finite
$d$-dimensional Minkowski content.  We establish this in
Section~\ref{sec:bounded-MC}.  We then transfer the half-plane
H\"older and occupation estimates through local conformal charts.
The Hausdorff- and packing-dimension identities require no angle
condition at the marked endpoints.  For the two box dimensions,
additional endpoint control is needed, and this is obtained under
the positive-angle assumption in Theorem~\ref{thm:smooth}.

As a further consequence of the H\"older and occupation estimates, we
also obtain the exact local dimension $d$ of the natural occupation
measure, uniformly over all points of the trace on the compact time
intervals covered by Proposition~\ref{prop:localdim}.

The remainder of the paper is organized as follows.
Section~\ref{sec:deterministic} develops the deterministic dimension criterion and reviews the
Minkowski-content and Minkowski-content-measure facts used below.
Section~\ref{sec:occupation} establishes the occupation estimates for
chordal and sssi SLE.
Section~\ref{sec:boundary} proves natural-time H\"older continuity for
chordal SLE in $\HH$ up to the boundary and concludes with the proofs
of the first two dimension theorems.
Section~\ref{sec:bounded-MC} establishes existence and finiteness of the Minkowski-content
parametrization for chordal SLE in bounded domains.
Section~\ref{sec:smooth} proves the bounded-domain dimension theorem
under the stated boundary assumptions.
Finally, Section~\ref{sec:localdim} records the local dimension of
occupation.

\section{Deterministic preliminaries} \label{sec:deterministic}
\subsection{Dimensions and a fixed-mesh criterion}
This section formulates a deterministic criterion for dimension
transformation using H\"older continuity and occupation bounds.
The proof uses the classical inverse-image covering argument underlying
Kaufman's dimension-doubling theorem; see, for example,
\cite[Section~9.4]{MortersPeres}.
In this subsection, $d$ is a deterministic parameter,
not necessarily related to $\kappa$.

We first recall the four notions of dimension used below;
see \cite[Chapters 2--3]{Falconer}.
For $E\subset\R^m$ and $s>0$, define
\BGE
 \mathcal H^s(E)
 =
 \lim_{\delta\downarrow0}
 \inf\left\{
 \sum_{j=1}^\infty(\diam U_j)^s:
 E\subset\bigcup_{j=1}^\infty U_j,\
 \diam U_j\le\delta
 \right\}. \label{Hausdorff-measure}
\EDE
Then $\mathcal H^s(E)$ is decreasing in $s$, and the Hausdorff dimension of $E$ is
\BGE
 \dimH E=\inf\{s>0:\mathcal H^s(E)=0\}. \label{Hausdorff-dim}
\EDE

For a nonempty bounded set $E$, let $\Nc(E,r)$ be the least
number of sets of diameter at most $r$ needed to cover $E$.
Its lower and upper box dimensions are
\BGE
 \dimBl E=\liminf_{r\downarrow0}
 \frac{\log\Nc(E,r)}{\log(1/r)},
 \qquad
 \dimBu E=\limsup_{r\downarrow0}
 \frac{\log\Nc(E,r)}{\log(1/r)}. \label{box-dim}
\EDE
These are also called lower and upper Minkowski dimensions.
When they agree, their common value is the Minkowski dimension.

We define the packing dimension of $E$ by
\begin{equation}\label{eq:packingchar}
 \dimP E
 =
 \inf\left\{
 \sup_{j\ge1}\dimBu E_j:
 E\subset\bigcup_{j=1}^\infty E_j,\
 E_j\subset\R^m\text{ bounded}
 \right\}.
\end{equation}
By \cite[(3.19) and Proposition~3.8]{Falconer}, this is equivalent
to the definition via packing measure. For the empty set, we set all four dimensions equal to zero.

Throughout the paper, we write $\L^1$ for Lebesgue measure on $\R$.
For a continuous map $\xi:J\to\R^m$ and a Borel set $I\subset J$, we write
\BGE \mu_I^\xi:=\xi_*(\ind_I\L^1) \label{mu}\EDE  for its occupation measure
on $I$.

Let $I$ be a nondegenerate compact real interval, and $f:I\to\R^m$ be continuous.

\begin{lemma}\label{lem:mesh}
Assume that, for some $0<\alpha\le1$ and $H\in(0,\infty)$,
\begin{equation}\label{eq:holderdet}
 |f(t)-f(s)|\le H|t-s|^\alpha,
 \qquad s,t\in I.
\end{equation}
Let $r>0$, and
\(
 h=(\frac{r}{2H})^{1/\alpha}
\).
 Partition $I$ into intervals of length $h$, allowing at most
one interval to be shorter. Let $z\in\R^m$, and let $N_{z,r}$ denote the number of mesh intervals which meet
$f^{-1}(\overline B(z,r))$. Then we have
\begin{equation}\label{eq:meshcount}
 N_{z,r}\le 1+h^{-1}\mu^f_I(B(z,2r)).
\end{equation}
\end{lemma}

\begin{proof}
If $J$ is one mesh interval that meets $f^{-1}(\overline B(z,r))$, i.e., there is $t\in J$ with $|f(t)-z|\le r$, then
\[
 |f(s)-z|\le |f(s)-f(t)|+|f(t)-z|\le Hh^\alpha+r=\frac32r<2r,\quad s\in J.
\]
Thus, $J$ is mapped into $B(z,2r)$. At least $N_{z,r}-1$ of the selected intervals have length $h$.
Since their interiors are disjoint,
$(N_{z,r}-1)h\le \mu^f_I(B(z,2r))$. This proves \eqref{eq:meshcount}.
\end{proof}

\begin{lemma} \label{prop:finite}
Assume \eqref{eq:holderdet}.  Suppose also that, for some $\beta\in (0, \alpha^{-1}]$ and $C,r_0>0$,
\begin{equation}\label{eq:massdet}
 \mu^f_I(B(z,r))\le Cr^\beta,
 \qquad z\in\R^m,\quad 0< r\le r_0.
\end{equation}
Let $\theta=\alpha^{-1}-\beta\ge 0$.
Then, for every $A\subset I$ and
$\mathfrak D \in\{\dimH,\dimP,\dimBl,\dimBu\}$,
\begin{equation}\label{eq:finitebound}
  \frac{ \mathfrak D (A)}{\alpha}-\theta
 \le \mathfrak D (f(A))\le \frac{\mathfrak D (A)}{\alpha}.
\end{equation}
\end{lemma}

\begin{proof}
The assertion is immediate for $A=\varnothing$, so assume
$A\ne\varnothing$.
The upper bounds for Hausdorff and box dimensions follow from the
standard H\"older covering estimates; see
\cite[Proposition~2.3 and Exercise~3.1]{Falconer}.
The packing-dimension bound then follows from
\eqref{eq:packingchar}.

By \eqref{eq:meshcount} and \eqref{eq:massdet}, for
$0<r\le(r_0/2)\wedge1$, the preimage of any closed ball
of radius $r$ is covered by at most $1+2^\beta C(2H)^{1/\alpha} r^{-\theta}$
intervals of length at most
$h:=(r/(2H))^{1/\alpha}$. Let
$C_0=1+2^\beta C(2H)^{1/\alpha}$. Then the number of intervals is no more than $C_0 r^{-\theta}$.
Every nonempty set of diameter at most $r$ is contained
in a closed ball of radius $r$. Consequently, $\Nc(A,h) \le C_0r^{-\theta}\Nc(f(A),r)$.
Taking logarithms, dividing by $\log(1/r)$, and then taking
the liminf and limsup as $r\downarrow0$ gives by \eqref{box-dim}
\[
 \frac1\alpha \mathfrak D (A)\le \mathfrak D (f(A))+\theta,
 \qquad \mathfrak D \in\{\dimBl,\dimBu\}.
\]
This gives the desired lower bounds for both $\dimBu$ and $\dimBl$.

For packing dimension, let $f(A)\subset\bigcup_j E_j$
be a countable cover by bounded sets, and put
$A_j=A\cap f^{-1}(E_j)$.
The inequality just proved for upper box dimension gives
\[
 \dimBu A_j
 \le\alpha\bigl(\dimBu f(A_j)+\theta\bigr)
 \le\alpha\bigl(\dimBu E_j+\theta\bigr).
\]
Since $A=\bigcup_j A_j$, \eqref{eq:packingchar} yields $\dimP A\le\alpha\bigl(\sup_j\dimBu E_j+\theta\bigr)$.
Taking the infimum over all such covers proves the
lower bound for packing dimension.

For Hausdorff dimension, let $s>\alpha(\theta+\dimH f(A))$. Then we may find $q>\dimH f(A)$ such that $s>\alpha(\theta+q)$. Let $\eps>0$. By \eqref{Hausdorff-measure} and \eqref{Hausdorff-dim}, one may choose a cover of $f(A)$ by closed balls
$\lin B(z_i,r_i)$ with $0<r_i\le (r_0/2)\wedge 1\wedge \eps$ and with
$\sum_i r_i^q\le \eps$.  The preimage of each $\lin B(z_i,r_i)$  is further covered by at most $C_0 r_i^{-\theta}$ intervals of length no more than $(\frac{r_i}{2H})^{1/\alpha}\le (2H)^{-1/\alpha}\eps^{1/\alpha}$. These intervals together form a cover of $A$, whose maximal diameter tends to $0$ as $\eps\to 0$, and whose $s$-dimensional cost is at most
\[
 \sum_i C_0 r_i^{-\theta} \Big(\frac{r_i}{2H}\Big)^{s/\alpha}=C_0 (2H)^{-s/\alpha} \sum_i r_i^{s/\alpha-\theta}\le C_0 (2H)^{-s/\alpha} \sum_i r_i^q\le C_0 (2H)^{-s/\alpha} \eps.
\]
Letting $\eps\downarrow 0$ and using \eqref{Hausdorff-measure}, we get $\mathcal H^s(A)=0$. Since this holds for any $s>\alpha(\theta+\dimH f(A))$, by \eqref{Hausdorff-dim},  one has $\dimH A\le \alpha(\theta+\dimH f(A))$, which implies the lower bound for $\dimH$.
\end{proof}

\begin{proposition}\label{cor:det}
Let $d\in[1,m]$. Suppose that for every $0<\alpha<1/d$
the map $f:I\to\R^m$ is $\alpha$-H\"older, and that for
every $0<\beta<d$ there are  constants $C_\beta,r_\beta\in (0,\infty)$ such that  $\mu_I^f(B(z,r))\le C_\beta r^\beta$
for all $z\in\R^m$ and all   $0<r<r_\beta$.
Then,  for every $A\subset I$,  \eqref{eq:main} holds with $f$ in place of $\eta$.
\end{proposition}

\begin{proof}
This follows from \eqref{eq:finitebound} by sending  $\alpha\uparrow1/d$ and $\beta\uparrow d$.
\end{proof}

\subsection{Minkowski content and natural parametrization}

We now return to the convention
$d=1+\kappa/8\in(1,2)$.
For a set $S\subset\C$ and $h>0$, write
\[
 S^h=\{z\in\C:\dist(z,S)<h\},
\]
and let $\L^2$ denote planar Lebesgue measure.  The upper and lower
$d$-dimensional Minkowski contents of $S$ are defined by
\[
 \lin\Cont_d(S)
 :=
 \limsup_{h\downarrow0}
 h^{d-2}\L^2(S^h),
 \qquad
 \ulin\Cont_d(S)
 :=
 \liminf_{h\downarrow0}
 h^{d-2}\L^2(S^h).
\]
When they agree, the common value is called the
$d$-dimensional Minkowski content of $S$ and is denoted by
$\Cont_d(S)$.

Let $J$ be a real interval.  We say that a continuous curve
$\xi:J\to\C$ is parametrized by $d$-dimensional Minkowski content if
\begin{equation}\label{eq:intervalcontent}
 \Cont_d(\xi(I))=\L^1(I),
 \qquad
 \text{for every compact interval }I\subset J.
\end{equation}
Throughout the paper, both the chordal SLE$_\kappa$ curve $\gamma$
and the sssi SLE$_\kappa$ curve $\eta$ are parametrized in this
sense.  For chordal SLE, the existence of this parametrization was
proved in \cite[Theorem~1.1]{LawlerRezaei}; for the sssi curve, see
\cite{ZhanHolder}.

We will also use the measure-valued version of Minkowski content.
Identifying $\C$ with $\R^2$, we follow
\cite[Definition~6.3]{ZhanBoundaryMC}.  Let $\mathcal M$ denote the
class of subsets of $\R^2$ which can be expressed as the intersection
of an open set with a closed set.  If $S\in\mathcal M$, a Borel
measure $M$ on $\R^2$ is called a $d$-dimensional Minkowski content
measure of $S$ if $M(\R^2\setminus S)=0$,
and, for every compact $K\subset S$,
\begin{equation}\label{eq:MC-finite}
 M(K)<\infty,
\end{equation}
and
\begin{equation}\label{eq:MC-continuity}
 M(\partial_S K)=0
 \quad\Longrightarrow\quad
 \Cont_d(K)=M(K).
\end{equation}
Here $\partial_S K$ denotes the relative boundary of $K$ in $S$.
Such a measure, when it exists, is unique by
\cite[Lemma~6.5]{ZhanBoundaryMC}; we denote it by
$\Md(\cdot;S)$.  Since every $S\in\mathcal M$ is $\sigma$-compact,
\eqref{eq:MC-finite} implies that $\Md(\cdot;S)$ is $\sigma$-finite.

We will repeatedly use the restriction property
\cite[Theorem~6.12]{ZhanBoundaryMC}.  If $\Md(\cdot;S)$ exists and
$G\subset\R^2$ is open, then
\begin{equation}\label{eq:MC-restriction}
 \Md(\cdot\cap G;S)
 =
 \Md(\cdot;S\cap G).
\end{equation}
We will also use the conformal covariance
\cite[Theorem~6.13]{ZhanBoundaryMC}: if $S$ is relatively closed in
an open set $G$, $\Md(\cdot;S)$ exists, and
$f:G\to\R^2$ is conformal and injective, then
\begin{equation}\label{eq:MC-covariance}
 f_*\bigl(|f'|^d\Md(\cdot;S)\bigr)
 =
 \Md(\cdot;f(S)).
\end{equation}

For chordal SLE $\gamma$, by \cite[Example~6.9]{ZhanBoundaryMC},  $\Md(\cdot;\gamma([0,\infty)))=\mu^\gamma_{[0,\infty)}$.
In \cite[Section~3.2]{LawlerRezaei}, this measure is called the
\emph{(natural) occupation measure} of the SLE curve. Likewise, for the sssi SLE $\eta$,
$\Md(\cdot;\eta(\R))=\mu^\eta_{\R}$;
see \cite[Section~3]{ZhanHolder}.

\section{Occupation estimates}\label{sec:occupation}

We first establish the occupation estimates used for the lower
dimension bounds.  We begin with chordal SLE in $\HH$, where the
required bounds follow from multipoint approach estimates.  We then
derive a corresponding multipoint estimate for sssi SLE and use it
to obtain the associated occupation bounds.  The complementary H\"older
estimate used for the upper dimension bounds will be proved in the
next section.

\subsection{Occupation estimates for chordal SLE}
\label{subsec:chordal-occupation}

Let $\R_+:=[0,\infty)$, and let
$\gamma:\R_+\to\overline{\HH}$ be chordal SLE$_\kappa$ in $\HH$
from $0$ to $\infty$, parametrized by $d$-dimensional Minkowski
content, and write $\Gamma=\gamma(\R_+)$.  We first obtain an
all-order local moment estimate for its natural occupation measure
and then deduce the almost-sure power bounds needed below.
Related higher-moment estimates on bounded sets were proved in
\cite{RezaeiZhan}; the estimate below is uniform in $z$.

\begin{lemma} \label{prop:tubemoment}
For every $p\in\N$, $z\in\C$, and $r>0$,
\begin{equation}\label{eq:tubemoment-mu}
 \EE\left[\mu^\gamma_{\R_+}(B(z,r))^p\right]
 \lesssim_p r^{dp}.
\end{equation}
\end{lemma}

\begin{proof}
Write $\mu:=\mu^\gamma_{\R_+}=\Md(\cdot;\Gamma)$.
Fix $z\in\C$. Since almost surely $\mu$ is $\sigma$-finite, by Fubini's theorem there is a deterministic set
$G_z\subset(0,\infty)$ of full Lebesgue measure such that
\begin{equation}\label{eq:good-radius}
\text{for every }\rho\in G_z,\quad\text{almost surely } \mu(\partial B(z,\rho))=0.
\end{equation}

Fix $\rho\in G_z$ and set $K_\rho:=\Gamma\cap\overline{B(z,\rho)}$.
Since chordal SLE is transient, $K_\rho$ is compact.  Moreover,
$\partial_\Gamma K_\rho
 \subset \Gamma\cap\partial B(z,\rho)$,
so \eqref{eq:good-radius} and
\eqref{eq:MC-continuity}   give
\begin{equation}\label{eq:good-radius-content}
 \mu(K_\rho)=\Cont_d(K_\rho)
 \qquad\text{a.s.}
\end{equation}

We next estimate the moments of the right-hand side.  For distinct
$w_1,\ldots,w_p\in\HH$, set $w_0=0$ and $\ell_j=\min_{0\le i<j}|w_j-w_i|$.
By \cite[Theorem~1.1 and (1.5)]{RezaeiZhan},
\begin{equation}\label{eq:chordal-multipoint}
 \PP\{w_j\in\Gamma^h,\ 1\le j\le p\}
 \lesssim_p
 \prod_{j=1}^p\left(\frac{h}{\ell_j}\right)^{2-d}.
\end{equation}
Consequently, for $0<h<\rho$,
\begin{align}
 &\EE\left[
 \left(
 h^{d-2}\L^2\bigl(B(z,2\rho)\cap\HH\cap\Gamma^h\bigr)
 \right)^p
 \right] \notag\\
 &\qquad\lesssim_p
 \int_{B(z,2\rho)^p}
 \prod_{j=1}^p\ell_j^{d-2}
 \,\prod_{j=1}^p\L^2(dw_j) \notag\\
 &\qquad\le
 \int_{B(z,2\rho)^p}
 \prod_{j=1}^p
 \left(\sum_{i=0}^{j-1}|w_j-w_i|^{d-2}\right)
 \,\prod_{j=1}^p\L^2(dw_j)
 \lesssim_p \rho^{dp}.
 \label{eq:good-radius-prelimit}
\end{align}
Here we used $\sup_{a\in\C}
 \int_{B(z,2\rho)}|w-a|^{d-2}\,\L^2(dw)
 \lesssim \rho^d $
and integrated successively in $w_p,\ldots,w_1$.

Since $K_\rho\subset\overline{\HH}$, for $0<h<\rho$,
\[
 K_\rho^h
 \subset
 \bigl(B(z,\rho+h)\cap\HH\cap\Gamma^h\bigr)
 \cup
 \bigl(B(z,\rho+h)\cap\{-h<\Im w\le0\}\bigr).
\]
The second set has area $O(\rho h)$, and hence, since $d>1$,
\[
 h^{d-2}
 \L^2\bigl(B(z,\rho+h)\cap\{-h<\Im w\le0\}\bigr)
 \longrightarrow0.
\]
As $B(z,\rho+h)\subset B(z,2\rho)$ for $h<\rho$, it follows from
the existence of $\Cont_d(K_\rho)$ that
\[
 \Cont_d(K_\rho)
 \le
 \liminf_{h\downarrow0}
 h^{d-2}
 \L^2\bigl(B(z,2\rho)\cap\HH\cap\Gamma^h\bigr).
\]
Fatou's lemma and \eqref{eq:good-radius-prelimit} therefore imply
\begin{equation}\label{eq:good-radius-moment}
 \EE[\mu(K_\rho)^p]
 =
 \EE[\Cont_d(K_\rho)^p]
 \lesssim_p \rho^{dp},
 \qquad \rho\in G_z.
\end{equation}

Finally, fix $r>0$ and choose a  sequence
$(\rho_n)$ in $G_z$ with $\rho_n\downarrow r$.  Since $\mu$ is supported
on $\Gamma$, $\mu(B(z,r)) \le \mu(K_{\rho_n})$.
Thus \eqref{eq:good-radius-moment} gives $\EE[\mu(B(z,r))^p]
 \lesssim_p \rho_n^{dp}$.
Letting $n\to\infty$ proves \eqref{eq:tubemoment-mu}.
\end{proof}

\begin{proposition} \label{prop:chordalmass}
There is a probability-one event, on which for every $R\in(0,\infty)$ and every
$0<\beta<d$, there is a random constant  $r_{R,\beta}(\gamma)\in(0,1]$ such that
\begin{equation}\label{eq:chordalmass}
 \mu^\gamma_{\R_+}(B(z,r))
 \le 2^d r^\beta,
 \qquad |z|\le R,\quad 0<r<r_{R,\beta}(\gamma).
\end{equation}
On the same event, for every $N>0$ and every $0<\beta<d$,
there is a finite random constant $C_{N,\beta}(\gamma)$ such that
\BGE
 \mu^\gamma_{[0,N]}(B(z,r))
 \le C_{N,\beta}(\gamma)r^\beta,
 \qquad z\in\C,\quad 0<r\le1. \label{eq:chordalmass2}
\EDE
\end{proposition}

\begin{proof}
Fix $R>0$ and $\beta\in(0,d)$.  At scale $2^{-n}$, there are
$O_R(2^{2n})$ lattice points of $2^{-n}(\Z+i\Z)$ in
$B(0,R+1)$.  Choose $p\in\N$ such that $ p(d-\beta)>2$.
By \eqref{eq:tubemoment-mu} and Markov's inequality,
\[
 \PP\left\{
 \mu^\gamma_{\R_+}(B(w,2^{2-n}))>2^{-n\beta}
 \text{ for some }
 w\in 2^{-n}(\Z+i\Z)\cap B(0,R+1)
 \right\}
 \lesssim_{R,p} 2^{-n(p(d-\beta)-2)}.
\]
The right-hand side is summable in $n$.  By the Borel--Cantelli
lemma, almost surely there is a finite random integer $n_0(\gamma)$ such that, for every
$n\ge n_0(\gamma)$ and every
$w\in 2^{-n}(\Z+i\Z)\cap B(0,R+1)$, $\mu^\gamma_{\R_+}(B(w,2^{2-n}))\le 2^{-n\beta}$.

Now let $|z|\le R$ and $0<r<2^{-n_0(\gamma)-1}$, and choose $n\ge n_0(\gamma)$
such that $ 2^{-n-1}<r\le 2^{-n}$.
If $w\in 2^{-n}(\Z+i\Z)$ is a nearest lattice point to $z$, then
$w\in B(0,R+1)$ and $B(z,r)\subset B(w,2^{2-n})$.
Consequently, $
 \mu^\gamma_{\R_+}(B(z,r))
 \le 2^{-n\beta}
 \le 2^d r^\beta$.
Thus \eqref{eq:chordalmass} holds with  $r_{R,\beta}=2^{-n_0(\gamma)-1}$ for any particular pair $(R,\beta)$.

Taking a countable intersection, we obtain a probability-one event
on which \eqref{eq:chordalmass} holds simultaneously for every
$R\in\N$ and every $\beta\in(0,d)\cap\Q$.  Work on this event.
For arbitrary $R>0$ and $0<\beta<d$, choose
$R'\in\N$ with $R'\ge R$ and
$\beta'\in(\beta,d)\cap\Q$.  Since
$r^{\beta'}\le r^\beta$ for $0<r\le1$, the estimate for
$(R',\beta')$ implies \eqref{eq:chordalmass} for $(R,\beta)$,
with $r_{R,\beta}(\gamma):=r_{R',\beta'}(\gamma)$.

Finally, fix $N>0$ and $0<\beta<d$, and choose an
integer-valued random variable $M_N(\gamma)\in\N$ such that every
ball of radius at most $1$ which meets $\gamma([0,N])$ has its
center in $B(0,M_N(\gamma))$. This is possible since $\gamma([0,N])$ is compact.
Set $r_0(\gamma):=r_{M_N(\gamma),\beta}(\gamma)$.
If $0<r<r_0(\gamma)$ and $B(z,r)$ meets $\gamma([0,N])$, then
$\mu^\gamma_{[0,N]}(B(z,r))
 \le
 \mu^\gamma_{\R_+}(B(z,r))
 \le2^d r^\beta$;
if the ball does not meet $\gamma([0,N])$, the left-hand side is
zero.  For $r_0(\gamma)\le r\le1$, $ \mu^\gamma_{[0,N]}(B(z,r))
 \le  N  \le  N r_0(\gamma)^{-\beta}r^\beta$.
Therefore \eqref{eq:chordalmass2} holds with
$C_{N,\beta}(\gamma)
 :=
 2^d\vee N r_0(\gamma)^{-\beta}$.
\end{proof}

As an immediate consequence, on the probability-one event in Proposition \ref{prop:chordalmass},
\begin{equation}\label{eq:zero-boundary-time}
 \mu^\gamma_{\R_+}(\R)=0.
\end{equation}
Indeed, fix $m\in\N$ and choose $\beta\in(1,d)$.  The interval $[-m,m]$ can be covered by $O_m(r^{-1})$
balls of radius $r\le 1$.  By
\eqref{eq:chordalmass}, $ \mu^\gamma_{\R_+}([-m,m])
 \lesssim_m  r^{\beta-1}$.
Letting $r\downarrow0$ gives
$\mu^\gamma_{\R_+}([-m,m])=0$.  Since
$\R=\bigcup_{m\in\N}[-m,m]$, \eqref{eq:zero-boundary-time} follows.

\subsection{Occupation estimates for sssi SLE}
\label{subsec:sssi-occupation}

Let $\eta:\R\to\C$ be an sssi SLE$_\kappa$, parametrized by
$d$-dimensional Minkowski content and normalized by $\eta(0)=0$.
It has self-similarity of index $1/d$ and stationary increments:
\begin{align}
 \bigl(\eta(ct)\bigr)_{t\in\R}
 &\overset{\mathrm{law}}{=}
 \bigl(c^{1/d}\eta(t)\bigr)_{t\in\R},
 &&c>0,\label{ss}\\
 \bigl(\eta(t+u)-\eta(u)\bigr)_{t\in\R}
 &\overset{\mathrm{law}}{=}
 \bigl(\eta(t)\bigr)_{t\in\R},
 &&u\in\R.\label{si}
\end{align}
Moreover, it satisfies reversibility: $(\eta(-t))_{t\in\R}\overset{\mathrm{law}}{=}(\eta(t))_{t\in\R}$, and is transient: $\lim_{t\to \pm \infty} \eta(t)=\infty$.
Write $\R_-=(-\infty,0]$ and $\Gamma_-:=\eta(\R_-)$.

We first derive for $\Gamma_-$ a multipoint approach estimate
analogous to the chordal estimate used in
Lemma~\ref{prop:tubemoment}.

\begin{proposition}\label{prop:sssi-multipoint}
For every $n\in\N$, let $z_1,\ldots,z_n\in\C\setminus\{0\}$ be
distinct, set $z_0=0$, and define
\[
 \lambda_j=\dist(z_j,\{z_0,\dots,z_{j-1}\}),
 \qquad 1\le j\le n.
\]
Then, for every $\eps>0$,
\begin{equation}\label{eq:sssi-multipoint}
 \PP\left\{
 \dist(z_j,\Gamma_-)<\eps,\ 1\le j\le n
 \right\}
 \lesssim_n
 \prod_{j=1}^n
 \left(1\wedge\frac{\eps}{\lambda_j}\right)^{2-d}.
\end{equation}
\end{proposition}

\begin{proof}
By the scaling relation \eqref{ss}, we may  assume that $ \max_{1\le j\le n}|z_j|<1$.
In particular, $\lambda_j\le1$ for every $j$.  If $\eps\ge1$,
the right side of \eqref{eq:sssi-multipoint} is bounded below by a
positive constant, so the assertion is immediate.  We henceforth
assume $0<\eps<1$.

Modulo parametrization, $\eta|_{\R_-}$ is a whole-plane
SLE$_\kappa(2)$ from $\infty$ to $0$; see
\cite[Section~2.1 and the beginning of Section~3]{ZhanHolder}.
Run this curve forward from $\infty$ toward $0$, and let $D_t$ denote
the component of the complement of the explored hull that contains
$0$.  Let $\tau$ be the first time such that $ \operatorname{crad}_{D_\tau}(0)=48$.

Write $D=D_\tau$, let $a$ be the tip at time $\tau$, and let $b$ be
the prime end of $D$ inherited from the original starting point
$\infty$.  By the domain Markov property of whole-plane SLE$_\kappa(2)$,
conditionally on the curve up to time $\tau$, the remaining curve is
a radial SLE$_\kappa(2)$ in $D$ from $a$ to $0$, with force point
$b$.

Let $f:(D;a,b)\longrightarrow(\HH;0,\infty)$ be conformal such that $|f'(0)|=1$.  Put $z_0=0$, $w_j=f(z_j)$, and $y_j=\Im w_j$, $0\le j\le n$. By conformal covariance of the
conformal radius, $48=\operatorname{crad}_D(0)
  = {\operatorname{crad}_{\HH}(w_0)}/|f'(0)| =2y_0$. Thus, $y_0=24$.
The Koebe quarter and distortion theorems give
\begin{equation}\label{eq:sssi-local-distortion}
B(0,12)\subset D,\quad f(B(0,2))
 \subset B(w_0,3)
\end{equation}
and
\begin{equation}\label{eq:sssi-local-bilip}
|z-z'|/2\le |f(z)-f(z')|\le 3 |z-z'|,
 \qquad z,z'\in B(0,2).
\end{equation}

Since $z_j\in B(0,1)$ and $\eps<1$, each
$B(z_j,\eps)$ is contained in $B(0,2)$ and is disjoint from the
explored hull.  Thus the event
$\{\dist(z_j,\Gamma_-)<\eps,\ 1\le j\le n\}$ is determined by
the remaining curve.  Moreover,
\eqref{eq:sssi-local-bilip} implies
\[
 f(B(z_j,\eps))
 \subset B(w_j,3\eps),
 \qquad 1\le j\le n.
\]

Let $\gamma$ be ordinary chordal SLE$_\kappa$ in $\HH$ from $0$ to
$\infty$, and write $H_t$ for its Loewner domain.  As long as
$w_0\in H_t$, set $ \Upsilon_t=\frac12\operatorname{crad}_{H_t}(w_0)$, and $\rho_\delta=\inf\{t:\Upsilon_t=\delta\}$. 
Thus $\Upsilon_0=y_0=24$.  Fix $\delta_2\in(0,1)$.  For
$0<\delta_1<3\delta_2/4$, condition $\gamma$ on
$\{\rho_{\delta_1}<\infty\}$ and stop it at $\rho_{\delta_2}$.
By \cite[Proposition~2.13]{LawlerWerness}, the law of this stopped
curve and that of two-sided radial SLE$_\kappa$ through $w_0$,
stopped at $\rho_{\delta_2}$, are mutually absolutely continuous,
and the Radon--Nikodym derivative of the latter with respect to the
former is
\[
 1+O\bigl((\delta_1/\delta_2)^u\bigr)
\]
uniformly, for some $u>0$.  The part of the two-sided radial curve
up to $w_0$ is radial SLE$_\kappa(2)$ in $\HH$ from $0$ to $w_0$,
with force point at $\infty$.  Hence, for fixed $\delta_2$, the
conditioned chordal law up to $\rho_{\delta_2}$ converges in total
variation, as $\delta_1\downarrow0$, to the corresponding radial
SLE$_\kappa(2)$ law.

Let $A_{\delta_2}$ be the event that this radial SLE$_\kappa(2)$,
before $\rho_{\delta_2}$, comes within distance $3\eps$ of each of
$w_1,\ldots,w_n$.  The total variation convergence gives
\begin{equation}\label{eq:sssi-conditioned-multipoint}
 \PP[A_{\delta_2}]
 \le
 \limsup_{\delta_1\downarrow0}
 \PP\left[
 \dist(\Gamma,w_j)<3\eps,\ 1\le j\le n
 \,\middle|\,
 \rho_{\delta_1}<\infty
 \right].
\end{equation}
As $\delta_2\downarrow0$, the stopping times $\rho_{\delta_2}$
increase to the terminal time at $w_0$.  Hence the events
$A_{\delta_2}$ increase to the event that the whole radial
SLE$_\kappa(2)$ comes within distance $3\eps$ of every $w_j$.
Therefore its probability is bounded by the right-hand side of
\eqref{eq:sssi-conditioned-multipoint}.

We use the notation of \cite{RezaeiZhan}.  Let
$a_\partial=8/\kappa-1$, and
\[
 P_y(x)=
 \begin{cases}
 y^{a_\partial-(2-d)}x^{2-d},&0\le x\le y,\\
 x^{a_\partial},&x\ge y.
 \end{cases}
\]
Set $L_0=|w_0|$, and $ L_j=\dist(w_j,\{0,w_0,\ldots,w_{j-1}\})$, $1\le j\le n$.  
By the Koebe quarter theorem,
\[
 \frac12\dist(w_0,\partial H_t)
 \le \Upsilon_t
 \le 2\dist(w_0,\partial H_t).
\]
Since $y_0=24$ and $\delta_1<1$, on
$\{\rho_{\delta_1}<\infty\}$ the nearest boundary point to $w_0$
cannot lie on $\R$.  Consequently, $\{\rho_{\delta_1}<\infty\}
 \subset \{\dist(\Gamma,w_0)\le 2\delta_1\}$.
The full multipoint estimate
\cite[Theorem~1.1]{RezaeiZhan} therefore gives
\begin{align*}
 &\PP\left[
 \rho_{\delta_1}<\infty,\,
 \dist(\Gamma,w_j)<3\eps,\ 1\le j\le n
 \right]\\
 &\qquad\lesssim_n
 \frac{P_{y_0}((2\delta_1)\wedge L_0)}
      {P_{y_0}(L_0)}
 \prod_{j=1}^n
 \frac{P_{y_j}((3\eps)\wedge L_j)}
      {P_{y_j}(L_j)}.
\end{align*}

On the other hand, \cite[Lemma 2.10]{LawlerWerness} gives, for $\delta_1$ small enough,
\begin{align*}
 \PP\{\rho_{\delta_1}<\infty\}
 & =
 c_*
 \left(\frac{\delta_1}{y_0}\right)^{2-d}
 \left(\frac{y_0}{L_0}\right)^{a_\partial}
 \left[1+O\left((\delta_1/y_0)^u\right)\right]\\
 &\asymp \left(\frac{2\delta_1}{y_0}\right)^{2-d}
 \left(\frac{y_0}{L_0}\right)^{a_\partial} = \frac{P_{y_0}((2\delta_1)\wedge L_0)}{P_{y_0}(L_0)},
\end{align*}
where the last equality follows because $2\delta_1<y_0\le L_0$. 
Thus, after conditioning on $\{\rho_{\delta_1}<\infty\}$, the factor
corresponding to $w_0$ cancels up to a constant independent of
$\delta_1$.

We next compare $L_j$ with $\lambda_j$.  By
\eqref{eq:sssi-local-bilip}, $\lambda_j/2\le \dist(w_j,\{w_0,\dots,w_{j-1}\})\le 3 \lambda_j\le 3$.  On the other hand,
\eqref{eq:sssi-local-distortion} and \eqref{eq:sssi-local-bilip} give $|w_j|\ge |w_0|-|w_j-w_0|\ge y_0-3= 21$. Thus,
\[
L_j=|w_j|\wedge \dist(w_j,\{w_0,\dots,w_{j-1}\})=\dist(w_j,\{w_0,\dots,w_{j-1}\})\ge \lambda_j/2.
\]

By \cite[(1.5)]{RezaeiZhan},
\[
 \frac{P_{y_j}((3\eps)\wedge L_j)}
      {P_{y_j}(L_j)}
 \le
 \left(1\wedge\frac{3\eps}{L_j}\right)^{2-d}
 \le
 6^{2-d}
 \left(1\wedge\frac{\eps}{\lambda_j}\right)^{2-d}.
\]
Hence the right-hand side of
\eqref{eq:sssi-conditioned-multipoint} is  $ \lesssim_n
 \prod_{j=1}^n
 \left(1\wedge\frac{\eps}{\lambda_j}\right)^{2-d}$.
Since the explored part of $\Gamma_-$ is disjoint from
$B(z_j,\eps)$ and $ f(B(z_j,\eps))\subset B(w_j,3\eps)$,
the same bound holds   for the original curve:
\[
 \PP\left[
 \dist(z_j,\Gamma_-)<\eps,\ 1\le j\le n
 \,\middle|\,\mathcal F_\tau
 \right]
 \lesssim_n
 \prod_{j=1}^n
 \left(1\wedge\frac{\eps}{\lambda_j}\right)^{2-d}.
\]
Averaging over $\mathcal F_\tau$ proves
\eqref{eq:sssi-multipoint}.
\end{proof}

\begin{remark}
The upper bound in \eqref{eq:sssi-multipoint} has exactly the same product
form as the corresponding multipoint estimate for ordinary
whole-plane SLE$_\kappa$; see \cite[Theorem~2]{MackeyZhan}.
However, $\Gamma_-$ is a whole-plane SLE$_\kappa(2)$, and hence is
a different process from ordinary whole-plane SLE$_\kappa$.
\end{remark}

We now use Proposition~\ref{prop:sssi-multipoint} exactly as the
chordal multipoint estimate was used in
Lemma~\ref{prop:tubemoment}.

\begin{lemma}\label{prop:sssi-tubemoment}
For every $p\in\N$, $z\in\C$, and $r>0$, $\EE [
 \mu^\eta_{\R}(B(z,r))^p] \lesssim_p r^{dp}$.
\end{lemma}

\begin{proof}
We first prove
\begin{equation}\label{eq:sssi-negative-moment}
 \EE\left[
 \mu^\eta_{\R_-}(B(z,r))^p
 \right]
 \lesssim_p r^{dp}.
\end{equation}
Write  $\mu_-:=\mu^\eta_{\R_-}=\Md(\cdot;\Gamma_-)$.  Fix $z\in\C$.  As in the proof of
Lemma~\ref{prop:tubemoment}, there is a deterministic  set $G_z\subset(0,\infty)$ of full Lebesgue measure such that,
\begin{equation}\label{eq:sssi-good-radius}
\text{for every }\rho\in G_z,\quad\text{almost surely } \mu_-(\partial B(z,\rho))=0.
\end{equation}

Fix $\rho\in G_z$ and set $K_\rho:=\Gamma_-\cap\overline{B(z,\rho)}$.
The transience of $\eta$ implies that $K_\rho$ is compact.
Moreover, $  \partial_{\Gamma_-}K_\rho
 \subset \Gamma_-\cap\partial B(z,\rho) $.
Hence, by \eqref{eq:sssi-good-radius} and
\eqref{eq:MC-continuity},
\begin{equation}\label{eq:sssi-good-radius-content}
 \mu_-(K_\rho)=\Cont_d(K_\rho)
 \qquad\text{a.s.}
\end{equation}

For distinct $w_1,\ldots,w_p\in\C\setminus\{0\}$, set $w_0=0$ and
$\lambda_j=\dist(w_j,\{w_0,\ldots,w_{j-1}\})$, $1\le j\le p$.
Proposition~\ref{prop:sssi-multipoint} gives
\[
 \PP\{w_j\in\Gamma_-^h,\ 1\le j\le p\}
 \lesssim_p
 \prod_{j=1}^p
 \left(1\wedge\frac{h}{\lambda_j}\right)^{2-d}\le \prod_{j=1}^p
 \left( \frac{h}{\lambda_j}\right)^{2-d}.
\]
Then the same argument as in \eqref{eq:good-radius-prelimit}  yields, for $0<h<\rho$,
\BGE
 \EE\left[
 \left(
 h^{d-2}\L^2(B(z,2\rho)\cap\Gamma_-^h)
 \right)^p
 \right]
 \lesssim_p \rho^{dp}.
 \label{eq:sssi-prelimitmoment}
\EDE
Since $ K_\rho^h
 \subset B(z,2\rho)\cap\Gamma_-^h$, $\Cont_d(K_\rho)
 \le
 \liminf_{h\downarrow0}
 h^{d-2}\L^2(B(z,2\rho)\cap\Gamma_-^h)$.
Combining this with
\eqref{eq:sssi-good-radius-content},
Fatou's lemma, and \eqref{eq:sssi-prelimitmoment}, we obtain $\EE[\mu_-(K_\rho)^p]
 \lesssim_p \rho^{dp}$ for each $\rho\in G_z$.

Finally, for arbitrary $r>0$, choose $\rho_n\in G_z$ with
$\rho_n\downarrow r$.  Since $\mu_-$ is supported on $\Gamma_-$,
 $\mu_-(B(z,r))\le \mu_-(K_{\rho_n})$,
and hence $\EE[\mu_-(B(z,r))^p]\lesssim_p \rho_n^{dp}$.
Letting $n\to\infty$ proves \eqref{eq:sssi-negative-moment}.

By reversibility of $\eta$,  \eqref{eq:sssi-negative-moment} holds with $\R_-$
replaced by $\R_+$. Since $\mu^\eta_{\R}  =
 \mu^\eta_{\R_-}+\mu^\eta_{\R_+}$,
the inequality
$(a+b)^p\le2^{p-1}(a^p+b^p)$ completes the proof.
\end{proof}

\begin{proposition}
\label{prop:sssi-occupation}
\label{prop:sssi-negative-mass}
There is a probability-one event on which, for every
$R\in(0,\infty)$ and every $0<\beta<d$, there is a random constant
$r _{R,\beta}(\eta)\in(0,1]$ such that
\begin{equation}\label{eq:sssi-negative-mass}
 \mu^\eta_{\R}(B(z,r))
 \le 2^d r^\beta,
 \qquad
 |z|\le R,\quad
 0<r<r _{R,\beta}(\eta).
\end{equation}
On the same event, for every  $N>0$ and
every $0<\beta<d$, there is a finite random constant
$C_{N,\beta}(\eta)$ such that
\begin{equation}\label{eq:sssi-occupation}
 \mu^\eta_{[-N,N]}(B(z,r))
 \le C_{N,\beta}(\eta)r^\beta,
 \qquad
 z\in\C,\quad 0<r\le1.
\end{equation}
\end{proposition}
\begin{proof}
The proof is the same as that of
Proposition~\ref{prop:chordalmass}, using
Lemma~\ref{prop:sssi-tubemoment}.
\end{proof}

\section{Natural-time H\"older continuity for chordal SLE}
\label{sec:boundary}

The upper bounds in the dimension theorems are obtained from
H\"older continuity.  For the sssi SLE$_\kappa$ process, the required
natural-time H\"older continuity is already known.  The main purpose
of this section is therefore to establish the corresponding result
for chordal SLE$_\kappa$ in $\HH$ from $0$ to $\infty$, parametrized
by $d$-dimensional Minkowski content.  At the end of the section, we
combine these H\"older estimates with the occupation estimates of
Section~\ref{sec:occupation} and the deterministic criterion of
Section~\ref{sec:deterministic} to prove the first two dimension
theorems.

 We continue to use the sssi curve $\eta$ and the chordal curve
$\gamma$ from the previous section. We couple them as follows. Write $\R_-=(-\infty,0]$.
Let  $U$ be the unbounded connected component of
$\C\setminus\eta(\R_-)$.
Conditionally on $\eta|_{\R_-}$, the curve $\eta|_{\R_+}$,
viewed modulo time parametrization, is chordal SLE$_\kappa$
in $U$ from $0$ to $\infty$; see
\cite[Section 2.1 and the beginning of Section 3]{ZhanHolder}.

Choose a conformal map $\varphi:\HH\to U$ sending $0$
and $\infty$ to the corresponding marked prime ends. It extends continuously to $\lin\HH$ in the spherical topology.
To fix its scale, choose the first point $z_*\in U$
in a fixed enumeration of $\Q+i\Q$, and require
$\Im\varphi^{-1}(z_*)=1$.
Thus $\varphi$ is determined by $\eta|_{\R_-}$.

We may therefore couple $\gamma$ and $\eta$
so that $\gamma$ is independent of $\eta|_{\R_-}$,
and hence of $\varphi$, and
\BGE
 \eta(\Theta(t))=\varphi(\gamma(t)),
 \qquad t\ge0, \label{Theta}
\EDE
where
\(
 \Theta(t):=\Cont_d\bigl(\varphi(\gamma([0,t]))\bigr)
\), $t\ge0$, is an increasing homeomorphism of
$\R_+$ onto itself.

For every $0\le u<v$ such that $\gamma([u,v])\subset\HH$,
the additivity of the natural occupation measure together with the
conformal covariance \eqref{eq:MC-covariance} gives
\begin{equation}\label{eq:coupled-clock}
 \Theta(v)-\Theta(u)
 =\Cont_d\bigl(\varphi(\gamma([u,v]))\bigr)
 =\int_u^v|\varphi'(\gamma(t))|^d\,\dd t.
\end{equation}

We record the known sssi H\"older estimate separately, since it will
be used directly in the proof of the first dimension theorem.  The
following is \cite[Theorem~1.2]{ZhanHolder}.

\begin{proposition}\label{prop:sssi-holder}
Almost surely, the sssi SLE$_\kappa$ process $\eta$ is locally
H\"older continuous of every order $\alpha\in(0,1/d)$.
\end{proposition}

A substantially sharper result was obtained recently by Holden and
Yuan \cite[Theorem~1.2]{HoldenYuan}, who established a critical
modulus-of-continuity estimate at the exponent $1/d$ with a
logarithmic correction.

Throughout the remainder of the paper, we work on the
probability-one event on which $\eta$ is transient:
$|\eta(t)|\to\infty$ as $t\to\pm\infty$,
and is locally H\"older continuous of every order
$\alpha\in(0,1/d)$, as in Proposition~\ref{prop:sssi-holder}.
In particular, $\eta(x+\R_-)$ is closed for every $x\in\R$.

For the proof of the chordal H\"older estimate, we need the following
moment estimate for $\eta$:
\begin{equation}\label{eq:wp-moments}
 \EE[\diam(\eta([0,1]))^q]
 \lesssim_q 1,\qquad q>0,
\end{equation}
which follows from \cite[Remark~1.8]{ZhanHolder}.

\subsection{Moving interior balls and independent conformal embedding}
Fix a constant $A\ge1$  such that the $\eta|_{\R_-}$-measurable event
\[
 G_A:=\{A^{-1}\le|\varphi'(i)|\le A,\ |\varphi(i)|\le A\}
 \quad\hbox{satisfies}\quad p_A:=\PP(G_A)>0.
\] For every event $E$ determined by
$\gamma$, independence gives $\PP(G_A\cap E)=p_A\PP(E)$.
Fix another constant $R\ge 2$. In the estimates below, constants may depend on $R$ and $A$, and such dependence will not be indicated. In particular, on $G_A$, $|\varphi'(i)|\asymp 1$ and $|\vphi(i)|\lesssim 1$.

\begin{lemma}\label{lem:ball-distortion}
There exist constants $c_0,C_0>0$ (depending on $R$ and $A$) such that,
for every $\rho\in (0,4R/3]$, and
$z=x+iy\in[-R,R]\times[3\rho/4,R]$, on $G_A$,
\begin{equation}\label{eq:polynomial-distortion}
 c_0 \rho< |\vphi'(z)|< C_0 \rho^{-3},
 \qquad
 \varphi(\overline B(z,\rho/4))\subset B(0,C_0 \rho^{-2}).
\end{equation}
Moreover, on $G_A$, for any $q,q'\in \overline B(z,\rho/4)$,
\begin{equation}\label{eq:local-ball-distortion}
 |\varphi'(q)|\asymp |\vphi'(z)|,\qquad
 |\varphi(q)-\varphi(q')|\gtrsim |\vphi'(z)||q-q'|.
\end{equation}
\end{lemma}

\begin{proof}
Put $\chi(\zeta)=i\frac{1+\zeta}{1-\zeta}$ and
$f(\zeta)=(\varphi(\chi(\zeta))-\varphi(i))/(2i\varphi'(i))$.
Then $f$ is univalent on the unit disk
$\D=\{\zeta\in\C:|\zeta|<1\}$ and satisfies
$f(0)=0$ and $f'(0)=1$.  Let $w=\chi^{-1}(z)\in\D$. Then
\[
 1-|w|^2=\frac{4y}{|z+i|^2},\qquad
|\vphi'(z)|=\frac{4|\varphi'(i)|}{|z+i|^2}|f'(w)|.
\]
Recall that we omit the dependence of constants on $R$ and $A$.
For $z\in [-R,R]\times [3\rho/4,R]$, we have $ |z+i|^2\asymp 1$ and $1+|w|\asymp 1$. Thus,
$1-|w|\asymp y$. The Koebe derivative bounds on $|f'(w)|$
give $y\lesssim |\vphi'(z)|/|\vphi'(i)|\lesssim y^{-3}$.   Since $|\vphi'(i)|\asymp 1$ on $G_A$ and $y\ge3\rho/4$, the first part of \eqref{eq:polynomial-distortion} follows with suitable
constants $c_0,C_0>0$. For $q\in \overline B(z,\rho/4)$, we have
$|\Re q|\le 4R/3$ and $\rho/2\le\Im q\le 4R/3$ since $\rho\le 4R/3$.
Applying the same argument with $w=\chi^{-1}(q)$ and using Koebe's upper growth bound on $|f(w)|$ gives $
 |\varphi(q)-\varphi(i)|/|\vphi'(i)|\lesssim \rho^{-2}$.
Since $|\varphi(i)|,|\vphi'(i)|\lesssim 1$ on $G_A$ and $\rho\lesssim 1$, we get the second part of \eqref{eq:polynomial-distortion}.

Since $B(z,y)\subset\HH$ and $\rho/4\le y/3$,  Koebe distortion provides the first part of \eqref{eq:local-ball-distortion}.
For the two-point lower bound, fix $q,q'\in \overline B(z,\rho/4)$.
Then $B(q',\rho/2)\subset\HH$ and $|q-q'|\le \rho/2$.
Koebe's lower growth bound, extended to $|q-q'|=\rho/2$
by continuity, gives
\(
 |\varphi(q)-\varphi(q')|
 \ge\tfrac14|\varphi'(q')|\,|q-q'|
\),
which together with the first part of \eqref{eq:local-ball-distortion} implies the second part.
\end{proof}

For $L,\rho,h>0$, let $\mathcal W(L,\rho,h)$ be the event
that there are   $0\le \til u<\til v<\infty$ with
\[
 |\eta(\til u)|< L,\qquad \til v-\til u\le h,\qquad
 |\eta(\til v)-\eta(\til u)|\ge\rho.
\]

\begin{lemma}\label{lem:wp-fast}
For every $M>0$,
\begin{equation}\label{eq:wp-fast}
 \PP[\mathcal W(L,\rho,h)]
 \lesssim_M (L/\rho)^d(h/\rho^d)^M,
 \qquad L\ge\rho>0,\quad h>0.
\end{equation}
\end{lemma}
\begin{proof}
We may assume that $h\le \rho^d$, since otherwise the right-hand
side of \eqref{eq:wp-fast} is at least $1$.
By the scaling relation \eqref{ss}, it suffices to consider $\rho=1$, so that
$L\ge1\ge h>0$.

For a path $f:\R\to\C$, set $D(f):=\diam(f([0,2h]))$.
By \eqref{ss} and \eqref{eq:wp-moments}, for every $q>0$,
\begin{equation}\label{eq:wp-Dmoment}
 \EE[D(\eta)^q]\lesssim_q h^{q/d}.
\end{equation}

Suppose   $(\til u,\til v)$ is a witness to
$\mathcal W(L,1,h)$.  If $\til u<h$, then
$\til u,\til v\in[0,2h]$, and hence $D(\eta)\ge1$.
Now suppose that $\til u\ge h$.  For $u\ge0$, write
$ (\vartheta_u\eta)(s):=\eta(u+s)-\eta(u)$, $s\in\R$.
For every $u\in[\til u-h,\til u]$, from   $\til u, \til v\in [u,u+2h]$, we get
$D (\vartheta_u\eta)\ge1$ and $ |\eta(u)|
 < L+D (\vartheta_u\eta) $.
Consequently,
\[
 Z(\eta):=\int_0^\infty
 \ind_{\{D (\vartheta_u\eta)\ge1\}}
 \ind_{\{|\eta(u)|< L+D (\vartheta_u\eta)\}}\L^1(du)\ge h.
\]
Thus, $\mathcal W(L,1,h)\subset\{D(\eta)\ge 1\}\cup \{Z(\eta)/h\ge 1\}$, which implies that
\begin{equation}\label{eq:wp-WZ}
 \PP[\mathcal W(L,1,h)]
 \le \PP(D(\eta)\ge1)+h^{-1}\EE Z(\eta).
\end{equation}

For each $u\ge0$, $ (\vartheta_u\eta)(-u)=-\eta(u)$,
and, by   \eqref{si},
$\vartheta_u\eta$ has the same law as $\eta$.  Hence, by Tonelli's
theorem, $\EE Z(\eta) =
 \EE\bigl[\ind_{\{D(\eta)\ge1\}}\mu^\eta_{\R_-}(B(0,L+D(\eta))\bigr]$.
Decomposing according to $2^{j-1}\le D(\eta)<2^{j}$, $j\in\N$, and using
Cauchy--Schwarz, $\mu^\eta_{\R_-}\le\mu^\eta_{\R}$,
Lemma~\ref{prop:sssi-tubemoment} with $p=2$ and $z=0$, and
\eqref{eq:wp-Dmoment}, we obtain, for $q>2d$,
$$ \EE Z(\eta)
 \le
 \sum_{j=1}^\infty
 \|\mu^\eta_{\R_-}(B(0,L+2^j))\|_{L^2}
 \PP(D(\eta)\ge2^{j-1})^{1/2}
 \lesssim_q
 L^d h^{q/(2d)}
 \sum_{j=1}^\infty 2^{j(d-q/2)}
 \lesssim_q L^d h^{q/(2d)}.
$$
Together with \eqref{eq:wp-WZ} and \eqref{eq:wp-Dmoment}, this gives
$ \PP[\mathcal W(L,1,h)]
 \lesssim_q
 h^{q/d}+L^d h^{q/(2d)-1}$.
Given $M>0$, choosing $q=2d(M+1)$  yields \eqref{eq:wp-fast} for $\rho=1$.
\end{proof}

For $z\in\C$ and $\rho,h>0$, let $\mathcal E(z,\rho,h)$ be the event that there exist $0\le u<v<\infty$ such that
\[
 v-u\le h,\quad \gamma([u,v])\subset \overline B(z,\rho /4),\quad
 |\gamma(v)-\gamma(u)|\ge \rho /8.
\]

\begin{lemma}\label{prop:ball-fast}
For  every $M,h >0$, $\rho\in (0,4R/3]$, and  $z\in [-R,R]\times [3\rho/4,R]$,
\begin{equation}\label{eq:ball-fast}
 \PP(\mathcal E(z,\rho,h))
 \lesssim_M \rho^{-4d}(h/\rho^d)^M.
 \end{equation}
\end{lemma}

\begin{proof}
Let $c_0,C_0$ be the constants in Lemma \ref{lem:ball-distortion}.
On $G_A\cap\mathcal E(z,\rho,h )$, let $(u,v)$ be a witness of $\mathcal E(z,\rho,h )$, and set
$\til u=\Theta(u)$ and $\til v=\Theta(v)$.
By \eqref{Theta},  \eqref{eq:coupled-clock}, and Lemma~\ref{lem:ball-distortion}, for some constants $C\ge 1\ge c>0$,
\[
\begin{aligned}
 \eta([\til u,\til v])
 &=\vphi(\gamma([u,v]))
 \subset B(0,C_0\rho^{-2}),\\
 \til v-\til u
 &\le C|\vphi'(z)|^dh ,\\
 |\eta(\til v)-\eta(\til u)|
 &\ge c|\vphi'(z)|\rho.
\end{aligned}
\]

 Let $\mathcal J
 =\left\{j\in\Z:
 [2^j,2^{j+1}]
 \cap[c_0\rho,C_0\rho^{-3}]\ne\varnothing
 \right\}$.
On $G_A$, by \eqref{eq:polynomial-distortion}, there is $j\in \mathcal J$ such that  $2^j\le|\vphi'(z)|\le2^{j+1}$. Then
$ \til v-\til u\le 4C2^{jd}h$ and $|\eta(\til v)-\eta(\til u)|
 \ge c 2^j\rho$.
Consequently,
\[
 G_A\cap\mathcal E(z,\rho,h )
 \subset
 \bigcup_{j\in\mathcal J}
 \mathcal W(C_0 \rho^{-2},c 2^j\rho,4C 2^{jd}h ).
\]

Set $\delta_0= c^d/(4C)<1$. For any  $j\in \mathcal J$, since $2^j\le C_0 \rho^{-3}$ and $c\le 1$, we get $C_0 \rho^{-2}\ge c 2^j \rho$. If $h /\rho^d\le\delta_0$, then
$4C 2^{jd}h \le(c2^j \rho)^d$, so
Lemma~\ref{lem:wp-fast} gives
\[
\PP\bigl(\mathcal W(C_0\rho^{-2},c2^j\rho,4C2^{jd}h )\bigr)\lesssim_M
 \left(\frac{C_0\rho^{-2}}{c 2^j\rho}\right)^d
 \left(\frac{4 C2^{jd}h }{(c 2^j\rho)^d}\right)^M\lesssim_M
 \rho^{-3d}2^{-jd}(h /\rho^d)^M.
\]
Since $2^j\ge c_0\rho/2$ for every $j\in\mathcal J$,
\(
 \sum_{j\in\mathcal J}2^{-jd}\lesssim \rho^{-d}
\).
A union bound and the independence of $G_A$ and
$\mathcal E(z,\rho,h )$ therefore imply
$ p_A\PP(\mathcal E(z,\rho,h ))
 =\PP(G_A\cap\mathcal E(z,\rho,h ))
 \lesssim_M \rho^{-4d}(h /\rho^d)^M$.
Dividing by the positive constant $p_A$  proves \eqref{eq:ball-fast} when
$h /\rho^d\le\delta_0$.
For $h/\rho^d>\delta_0$, \eqref{eq:ball-fast} trivially holds since its RHS is bounded below by a positive constant.
\end{proof}

\subsection{From crossing estimates to H\"older continuity}
Fix $R\ge 2$. We continue to allow constants to depend implicitly on $R$.
For $r,\rho>0$, let $\mathcal F(r,\rho)$
be the event that there exist $0\le u<v<\infty$ such that
\[
 \gamma([u,v])\subset[-R-2,R+2]\times[0,\rho],
 \qquad
 |\Re\gamma(v)-\Re\gamma(u)|\ge r.
\]

\begin{lemma}\label{lem:flat-strip}
Let $a_\partial=8/\kappa-1>0$. Then for any $m\in\N$,
\begin{equation}\label{eq:flat-strip}
 \PP(\mathcal F(r,\rho))
 \lesssim_m r^{-1}(\rho/r)^{a_\partial m}, \quad r,\rho>0.
\end{equation}
\end{lemma}
\begin{proof}
We may assume that $r\le 2(R+2)$ for otherwise $\mathcal F(r,\rho)$ cannot happen.
For $k\in\Z$, let $I_k=[kr/4,(k+1)r/4]$, and $E_k$ be the event that there exist  $0\le s<t$
such that
\[
 \gamma([s,t])\subset I_k\times[0,\rho],
 \qquad
 \{\Re\gamma(s),\Re\gamma(t)\}
 =\{kr/4,(k+1)r/4\}.
\]  Let $S$ be the set of $k\in\Z$ such that $I_k\subset[-R-2,R+2]$. Then $|S|\lesssim r^{-1}$.

Suppose that $\mathcal F(r,\rho)$ occurs, and  choose
$u<v$ witnessing this event. Let $J$ denote the interval between
$\Re\gamma(u)$ and $\Re\gamma(v)$. Then for some $k\in S$, $J$ contains  $I_k$. From $\gamma([u,v])\subset \R\times [0,\rho]$, we see that $\gamma([u,v])$ has some subarc that crosses $I_k\times [0,\rho]$ horizontally. Thus, $\mathcal F(r,\rho)\subset \bigcup_{k\in S} E_k$.

Fix $k$. We may pick $m$ deterministic points $x_1,\dots,x_m$ on $I_k$ such that  $\min_{0\le i<j\le m} |x_j-x_i|\ge r/(8m)$, where $x_0:=0$. On the event $E_k$, $\gamma$ visits all $\lin B(x_j,\rho)$, $1\le j\le m$. By \cite[Theorem 1.1]{RezaeiZhan},
$$\PP[E_k]\lesssim_{m} \prod_{j=1}^m \Big(\frac{\rho}{\dist(x_j,\{x_i:0\le i\le j-1\})}\Big)^{a_\pa}\lesssim_m (\rho/r)^{a_\pa m}.$$ Since $\mathcal F(r,\rho)\subset \bigcup_{k\in S} E_k$, and $|S|\lesssim r^{-1}$, a union bound gives \eqref{eq:flat-strip}.
\end{proof}

Let $\mathcal D_R(h,r)$ be the event that there exist $0\le s<t<\infty$ such that
\BGE
 |\gamma(s)|\le R,\qquad t-s\le h,\qquad
 |\gamma(t)-\gamma(s)|\ge r. \label{event-D}
\EDE

\begin{lemma}\label{lem:fast-dichotomy}
 For $0<\rho<r/2<1/2$, there exists a set $S\subset [-R-2,R+2]\times [3\rho/4,R+2]$ with $\# S\lesssim  \rho^{-2}$ such that
\begin{equation}\label{eq:fast-dichotomy}
 \mathcal D_R(h,r)
 \subset \mathcal F(r/2,\rho)
       \ \cup\ \bigcup_{z\in S}\mathcal E(z,\rho,h),\quad h>0.
\end{equation}
\end{lemma}
\begin{proof} Suppose that $\mathcal D_R(h,r)$ occurs, and choose
$s<t$ witnessing this event. By replacing $t$ with the
first time after $s$ at which the displacement from
$\gamma(s)$ equals $r$, we may assume that
$ |\gamma(t)-\gamma(s)|=r$ and
$ \gamma([s,t])\subset\overline B(\gamma(s),r)
 \subset\overline B(0,R+1)$.
If $\gamma([s,t])\subset\R\times[0,\rho]$, then
\(
 |\Re\gamma(t)-\Re\gamma(s)|
 \ge\sqrt{r^2-\rho^2}>r/2
\),
so $\mathcal F(r/2,\rho)$ occurs.
Otherwise, choose $v\in[s,t]$ such that
$\Im\gamma(v)>\rho$. At least one of $\gamma(s)$ and
$\gamma(t)$ is at distance at least $r/2$ from $\gamma(v)$.
Starting at $v$ and moving in time towards that endpoint,
let $u$ be the first time encountered at which
$|\gamma(u)-\gamma(v)|=\rho/8$.
The subsegment with time interval
$[u\wedge v, u\vee v]$ stays in
$\overline B(\gamma(v),\rho/8)$, has endpoint displacement
$\rho/8$, and has duration at most $h$.

A nearest point of $(\rho/64)\Z^2$ to $\gamma(v)$ is at
distance less than $\rho/64$. Since
$|\gamma(v)|\le R+1$ and $\Im\gamma(v)>\rho$,
this point, denoted by $z$, belongs to $S:=(\rho/64)\Z^2\cap ([-R-2,R+2]\times[3\rho/4,R+2])$. Note that $\# S\lesssim  \rho^{-2}$.
Then
\(
 \overline B(\gamma(v),\rho/8)
 \subset\overline B(z ,\rho/4)
\),
so $\mathcal E(z ,\rho,h)$ occurs.
This proves \eqref{eq:fast-dichotomy}.
\end{proof}

\begin{lemma} \label{thm:fast-bound}
 Let $\alpha\in (0,1/d)$, and $h_0=2^{-\frac 2{1/d-\alpha}}\in (0,1)$. Then for any $N>0$,
\begin{equation}\label{eq:fast-bound}
  {
 \PP(\mathcal D_R(h,h^\alpha))\lesssim_{\alpha,N} h^N,
 \qquad 0<h<h_0.}
\end{equation}
\end{lemma}
\begin{proof}
Let $\nu=(1/d+\alpha)/2$. Let $h\in (0,h_0)$. Set $r=h^\alpha$ and  $\rho=h^\nu$.  Then $\rho<r/2$ and $h\le \rho^d$.
Lemmas~\ref{prop:ball-fast}, \ref{lem:flat-strip}, and \ref{lem:fast-dichotomy},  with $R$ replaced by $R+2$
in the first lemma, give for any $m\in\N$ and $M>0$,
\[
 \PP(\mathcal D_R(h,h^\alpha))
 \lesssim_{m,M}  h^{a_\partial m(\nu-\alpha)-\alpha}
       + h^{M(1-\nu d)-\nu(4d+2)}.
\]
Since $a_\partial>0$, $\nu>\alpha$ and $\nu d<1$, we may choose $m,M\in\N$ to make  $a_\partial m(\nu-\alpha)-\alpha\ge N$ and $M(1-\nu d)-\nu(4d+2)\ge N$.  Then the proof is complete.
\end{proof}

\begin{proposition}\label{prop:chordal-holder}
Almost surely,  $\gamma$ is locally H\"older continuous of every order $\alpha\in (0,1/d)$.
\end{proposition}
\begin{proof}
Fix an integer $R\ge2$ and $0<\alpha<1/d$.
Applying \eqref{eq:fast-bound} with $N=1$ gives
\[
 \sum_{n=1}^{\infty}
 \PP\bigl(\mathcal D_R(2^{-n},2^{-n\alpha})\bigr)<\infty.
\]
By Borel--Cantelli, there is an almost surely finite
random integer $n_{R,\alpha}(\gamma)$ such that for any $n\ge n_{R,\alpha}(\gamma)$,
$\mathcal D_R(2^{-n},2^{-n\alpha})$ does not occur.
Taking a countable intersection over integers $R\ge2$
and a fixed sequence $0<\alpha_j\uparrow1/d$ gives
a common probability-one event on which this holds
for every $R\in\N_{\ge 2}$ and $\alpha_j$.

Suppose the above event occurs. Fix  a deterministic $T\in(0,\infty)$ and define random numbers $ M(\gamma)=\sup_{0\le t\le T}|\gamma(t)|<\infty$ and
$ R(\gamma)=\lceil M(\gamma)\rceil+3$. Define $\delta_{j}(\gamma)=2^{-n_{R(\gamma),\alpha_j}(\gamma)}$, $j\in\N$.
If $s,t\in[0,T]$ and $0<t-s\le\delta_{j}(\gamma)$,
choose $n\in\N$ such that $2^{-n-1}< t-s\le 2^{-n}$. Then $n=\lfloor -\log_2(t-s)\rfloor\ge   n_{R(\gamma),\alpha_j}(\gamma)$ since $-\log_2(t-s)\ge n_{R(\gamma),\alpha_j}(\gamma)$. The absence of
$\mathcal D_{R(\gamma)}(2^{-n},2^{-n \alpha_j})$ gives $ |\gamma(t)-\gamma(s)|
 <2^{-n\alpha_j}\le 2^{\alpha_j}(t-s)^{\alpha_j}$.
If $t-s>\delta_{j}(\gamma)$, then we use the bound
$ |\gamma(t)-\gamma(s)|\le 2M(\gamma)$.
Consequently, $\gamma$ is $\alpha_j$-H\"older continuous
on $[0,T]$ with  constant  $2^{\alpha_j}\vee(2M(\gamma)\delta_{j}(\gamma)^{-\alpha_j})$.
For any $0<\alpha<1/d$, choosing $j$ with
$\alpha<\alpha_j$ shows that $\gamma$ is
$\alpha$-H\"older continuous on $[0,T]$.
\end{proof}

\subsection{Proofs of the first two dimension theorems}\label{sec:proofmain}
\begin{proof}[Proof of Theorem~\ref{thm:main}]
Intersect the probability-one events provided by
Propositions~\ref{prop:sssi-holder} and
\ref{prop:sssi-occupation}.
On this common event, for every compact interval $I\subset\R$,
the required H\"older and occupation bounds hold for all
$\alpha\in(0,1/d)$ and $\beta\in(0,d)$.

For each integer $n\ge1$, Proposition~\ref{cor:det}
applies to $f=\eta|_{[-n,n]}$.
Consequently, simultaneously for every $A\subset[-n,n]$, \eqref{eq:main} holds.
Every bounded subset of $\R$ is contained in $[-n,n]$ for some
$n$, proving \eqref{eq:main} for bounded sets.

For arbitrary $A\subset\R$, countable stability of Hausdorff and
packing dimensions gives
\[
 \mathfrak D(\eta(A))
 =\sup_{n\ge1}\mathfrak D\bigl(\eta(A\cap[-n,n])\bigr)
 =d\sup_{n\ge1}\mathfrak D(A\cap[-n,n])
 =d\,\mathfrak D(A),
 \qquad
 \mathfrak D\in\{\dimH,\dimP\}.
\]
This proves the theorem.
\end{proof}

\begin{proof}[Proof of Theorem~\ref{thm:chordal}]
Intersect the probability-one events provided by
Propositions~\ref{prop:chordal-holder} and \ref{prop:chordalmass}.
On this common event, Proposition~\ref{cor:det} applies to
$f=\gamma|_{[0,n]}$ for every integer $n\ge1$,
giving all four dimension identities simultaneously
for every $A\subset[0,n]$.
The bounded-set conclusion and the extension to arbitrary
time sets follow exactly as in the proof of
Theorem~\ref{thm:main}.
\end{proof}

\section{Minkowski content parametrization in bounded domains}
\label{sec:bounded-MC}

We show in this section that the natural parametrization of chordal
SLE in the upper half-plane induces a parametrization by
$d$-dimensional Minkowski content after conformal transport to a
bounded domain with locally connected boundary, under a mild size assumption on its boundary.
Conformal covariance directly identifies the Minkowski content measure
on the part of the trace lying in the interior of the domain.  The
main point is to globalize this identification across the boundary.

\begin{lemma}\label{lem:MC-globalization}
Let $d\in(1,2)$, let $D\subset\C$ be a bounded domain, let $S\subset\overline D$ be compact such that $\Md(\cdot;S\cap D)$ exists.
Suppose that $\Cont_d(\pa D)=0$ and
\begin{equation}\label{eq:boundary-tightness}
 \lim_{r\downarrow0}\,
 \limsup_{h\downarrow0}
 h^{d-2}\L^2\bigl(S^h\cap(\partial D)^r\cap D\bigr)=0.
\end{equation}
Then $\Md(\cdot;S\cap D)=\Md(\cdot;S)$ and $\Md(S;S)<\infty$.
\end{lemma}

\begin{proof}
We write $M:=\Md(\cdot;S\cap D)$. Since $M$ is $\sigma$-finite, there are at most countably many $r>0$ such that $M(L_r)>0$, where   $L_r:= \{z\in S\cap D:\dist(z,\pa D)=r\}$.
Hence we may choose a sequence $r_n\downarrow0$ such that
\begin{equation}\label{eq:good-radii}
 M(L_{r_n})
 =
 M(L_{r_n/2})
 =0,\quad n\in\N.
\end{equation}

From $M(\C\sem (S\cap D))=0$, we get $M(\C\sem S)=0$. We need to verify  \eqref{eq:MC-finite} and \eqref{eq:MC-continuity}.
We first prove $M(S)<\infty$, which implies   \eqref{eq:MC-finite}.    Set
\[
 T(r):=
 \limsup_{h\downarrow0}
 h^{d-2}\L^2\bigl(S^h\cap(\partial D)^r\cap D\bigr).
\]
By \eqref{eq:boundary-tightness}, $T(r)\to0$ as $r\downarrow0$. Choose $m$ sufficiently large that $T(2r_m)<\infty$.  For $n>m$,
let
\[
 A_{r_n,r_m}
 :=
 \{z\in S\cap D:r_n\le\dist(z,\pa D)\le r_m\}.
\]
Then $A_{r_n,r_m}$ is a compact subset of $S\cap D$, and
$ \partial_{S\cap D}A_{r_n,r_m}\subset L_{r_n}\cup L_{r_m}$.
Hence, by \eqref{eq:good-radii}, $\Cont_d(A_{r_n,r_m})=M(A_{r_n,r_m})$.
Moreover, if $0<h<r_n$, then
$A_{r_n,r_m}^h\subset S^h\cap(\partial D)^{2r_m}\cap D$,
and therefore $M(A_{r_n,r_m})=\Cont_d(A_{r_n,r_m}) \le T(2r_m)$.
Letting $n\to\infty$   and using continuity from below gives
$ M  (\{z\in S\cap D:\dist(z,\pa D)\le r_m\})  \le T(2r_m)<\infty $.
On the other hand, $\{z\in S\cap D:\dist(z,\pa D)\ge r_m\}$
is a compact subset of $S\cap D$, and hence has finite $M$-mass.
Since $M(S\cap\partial D)=0$,
we conclude that $M(S)<\infty$.

We next verify \eqref{eq:MC-continuity}. Let $K\subset S$ be compact and suppose that $M(\partial_S K)=0$.
We need to prove $\Cont_d(K)=M(K)$. For $r>0$, set $K_r:=K\setminus(\partial D)^r$.
Then $K_r$ is a compact subset of $S\cap D$, and $\partial_{S\cap D}K_r\subset (\partial_S K\cap D) \cup L_r$.
Thus, for $r=r_n$ or $r=r_n/2$,  $M(\partial_{S\cap D}K_r)=0$
by  $M(\partial_S K)=0$ and \eqref{eq:good-radii}, which implies that
\begin{equation}\label{eq:Kr-content}
 \Cont_d(K_r)=M(K_r),\quad\text{for }r=r_n\text{ or }r_n/2.
\end{equation}
Thus, $\ulin\Cont_d(K)\ge\Cont_d(K_{r_n})=M(K_{r_n})$.
Sending $n\to\infty$, and using that  $K_{r_n}\uparrow K\cap D$ and $M(S\cap\partial D)=0$, we get  $\ulin\Cont_d(K)\ge M(K)$.

For the reverse inequality, fix $r_n$ and let $0<h<r_n/2$.
Then
$$K^h \subset K_{r_n/2}^{h}\cup\bigl(S^h\cap D\cap(\pa D)^{r_n}\bigr)\cup (\pa D)^h.$$
It follows that
$$
h^{d-2}\L^2(K^h)
 \le
h^{d-2}\L^2(K_{r_n/2}^{h})
 +
h^{d-2}\L^2\bigl(
S^h\cap D\cap(\pa D)^{r_n}
\bigr)
+
h^{d-2}\L^2((\pa D)^h).
$$
Taking $\limsup_{h\downarrow0}$ and using
\eqref{eq:Kr-content} and $\Cont_d(\pa D)=0$, we get
\[
\lin\Cont_d(K)
\le
M(K_{r_n/2})+T(r_n)
\le
M(K)+T(r_n).
\]
Letting $n\to\infty$ and using $T(r_n)\to 0$ yields $\lin\Cont_d(K)\le M(K)$.

Combining  $\ulin\Cont_d(K)\ge M(K)$ and $\lin\Cont_d(K)\le M(K)$, we obtain $\Cont_d(K)=M(K)$. Since this holds for any compact $K\subset S$ with $M(\partial_S K)=0$, we conclude that $M=\Md(\cdot;S)$. 
\end{proof}

We now apply this globalization lemma to conformal images of chordal
SLE. Recall that $\gamma$ is a chordal
SLE$_\kappa$ in $\HH$ from $0$ to $\infty$, parametrized by
$d$-dimensional Minkowski content, and $\Md(\cdot;\Gamma)= \mu^\gamma_{\R_+}$, where
$\Gamma:=\gamma(\R_+)$. Applying \eqref{eq:MC-restriction} to $G=\HH$ and using \eqref{eq:zero-boundary-time} yields
\begin{equation}\label{eq:MC-half-plane-restriction}
 \mu^\gamma_{\R_+} =\Md(\cdot;\Gamma\cap\HH) .
\end{equation}

Now let $D$ be a bounded simply connected domain with locally connected boundary and let
$F:\HH\to D$ be conformal. By the Carath\'eodory--Torhorst theorem
\cite[Theorem~2.1]{Pommerenke}, $F$ admits a continuous extension from the spherical closure
of $\HH$ onto $\overline D$.  We use the same notation for the extension and set  $a=F(0)$ and $b=F(\infty)$. We extend $\gamma$ to $[0,\infty]$ such that $\gamma(\infty)=\infty$. Set $\gamma_F:=F\circ\gamma$ and $\Gamma_D:=\gamma_F([0,\infty])$.
Then $\Gamma_D$ is a compact subset of $\overline D$. Since $\Gamma\cap\HH$ is relatively closed in $\HH$, and $F(\Gamma\cap\HH)=\Gamma_D\cap D$, by
\eqref{eq:MC-half-plane-restriction} and the conformal covariance
\eqref{eq:MC-covariance},
\begin{equation}\label{eq:target-local-MC}
 F_*\bigl( |F'|^d\,\mu^\gamma_{\R_+}
 \bigr)=\Md(\cdot;\Gamma_D\cap D).
\end{equation}

\begin{lemma}\label{lem:target-boundary-tight}
Suppose that, for some $q<d$, $\lin\Cont_q(\pa D)<\infty$. Then  $\Gamma_D$ a.s.\ satisfies
\eqref{eq:boundary-tightness}.
\end{lemma}

\begin{proof}
By the assumption, there exists a constant $r_0>0$ such
that
\begin{equation}\label{eq:boundary-q-tube}
 \L^2\bigl((\pa D)^r\bigr)
 \lesssim r^{2-q},
 \qquad 0<r\le r_0.
\end{equation}

We use the following dyadic Whitney decomposition $\mathcal W$ of
$\HH$.  For $n,m\in\Z$, let
\[
 Q_{n,m}
 :=
 [m2^{-n},(m+1)2^{-n}]
 \times[2^{-n},2^{1-n}],
\]
and let $\mathcal W=\{Q_{n,m}:n,m\in\Z\}$.  Thus every
$Q\in\mathcal W$ is a dyadic square with $ \dist(Q,\R)$ equal to its side length $\ell_Q$.
The interiors of the squares are pairwise disjoint, and the squares
cover $\HH$. Let $z_Q$ be the center of
$Q$, and put $a_Q:=|F'(z_Q)|$ and $\rho_Q:=\diam(F(Q))$.
By the Koebe distortion theorem, there is a constant $C_1>0$ such that, for any $Q\in\mathcal W$,
\BGE F(Q)\subset (\pa D)^{C_1\rho_Q},\quad \rho_Q\le C_1 \dist(F(Q),\pa D),\quad \L^2(F(Q)) \ge \rho_Q^2/C_1.\label{CQW}\EDE

For every $Q\in\mathcal W$, define $\widehat Q=z_Q+1.5(Q-z_Q)$. Then $\ha Q$ is the square
with the same center as $Q$, sides parallel to the coordinate axes,
and side length $1.5\ell(Q)$. Let $\mathcal N(Q)=\{R\in\mathcal W:R\cap \ha Q\ne \varnothing\}$. Then $\# \mathcal N(Q)=9$, and $\ell(R)\asymp \ell(Q)$ for every $R\in\mathcal N(Q)$.
For $R\in\mathcal N(Q)$,   write
\[
 Z_R
 :=
 \Cont_d^+(\Gamma\cap R;-\log(\ell(R))):=\sup_{s\ge -\log(\ell(R))} e^{(2-d)s}\L^2 ((\Gamma\cap R)^{\exp(-s)}),
\]
as in \cite[p.~1100]{LawlerRezaei}.  In particular,
\begin{equation}\label{eq:ZR-maximal}
 \eps^{d-2}
 \L^2\bigl((\Gamma\cap R)^\eps\bigr)
 \le Z_R,
 \qquad 0<\eps\le\ell(R).
\end{equation}
Let $G_{\HH}$ be the Green's function for $\gamma$ (cf.\ \cite[(8)]{LawlerRezaei}), which satisfies $G_{\HH}(z)\lesssim (\Im z)^{d-2}$.
By \cite[Lemma~3.5]{LawlerRezaei},
\begin{equation}\label{eq:ZR-first-moment}
 \EE(Z_R)\lesssim G(R):=\int_R G_{\HH}(z)\L^2(dz)\lesssim \int_R (\Im z)^{d-2}\L^2(dz)\asymp \ell(R)^d.
\end{equation}

Koebe distortion on $\widehat Q$ gives $|F'(z)|\asymp a_Q$ for $z\in\widehat Q$, and
$\dist (F(Q),\pa F(\widehat Q) )
 \gtrsim a_Q\ell_Q
 \asymp \rho_Q$.
Consequently, there is a constant $c_1>0$ such that, whenever
$0<h<c_1\rho_Q$ and $w\in F(Q)\cap\Gamma_D^h$,
a point of $\Gamma_D$ within distance $h$ of $w$ must belong to
$F(\Gamma\cap\widehat Q)$.  This implies by Koebe distortion again that, for some constant $C_2>0$,
 $F^{-1}(w) \in (\Gamma\cap \ha Q)^{C_2h/a_Q}$. Thus,
 $$ F(Q)\cap\Gamma_D^h\subset F((\Gamma\cap \ha Q)^{C_2h/a_Q})\subset \bigcup_{R\in\mathcal N(Q)}  F((\Gamma\cap R)^{C_2h/a_Q}).$$
For any $R\in\mathcal N(Q)$, since $\rho_Q\asymp a_Q \ell(Q)\asymp a_Q \ell(R)$, we can find a constant $c_2>0$ such that, if $h<c_2\rho_Q$, then
 $C_2h/a_Q\le\ell(R)/4$ for any $R\in\mathcal N(Q)$. Applying Koebe distortion first on $\ha Q$ and then on $\ha R:=z_R+1.5(R-z_R)$ for $R\in\mathcal N(Q)$, we get $|F'(z)|\asymp |F'(z_Q)|=a_Q$ for any $z\in (\Gamma\cap R)^{C_2h/a_Q}$.
Changing variables under $F$ and using \eqref{eq:ZR-maximal} gives,
\begin{align}
  h^{d-2}\L^2(F(Q)\cap\Gamma_D^h)&\lesssim a_Q^d \sum_{R\in\mathcal N(Q)}  (C_2h/a_Q)^{d-2}  \L^2 ((\Gamma\cap R)^{C_2h/a_Q}) \notag \\ &\le a_Q^d \sum_{R\in\mathcal N(Q)} Z_R=:Y_Q,\qquad \text{if }0<h<c_2\rho_Q.\label{YQ1}
\end{align}
 Since $\#\mathcal N(Q)=9$, by \eqref{eq:ZR-first-moment} and that $\ell(R)\asymp \ell(Q)$, $\EE(Y_Q)\lesssim a_Q^d \ell(Q)^d\asymp \rho_Q^d$.

For $j\in\Z$ and $h>0$, define
$$ \mathcal W_j :=\{Q\in\mathcal W:  2^{-j-1}<\rho_Q\le 2^{-j}\}, $$
$$\mathcal W^>_h :=\{Q\in\mathcal W:\rho_Q>h/c_2\},\quad \mathcal W^\le_h :=\{Q\in\mathcal W:\rho_Q\le h/c_2\},$$
where $c_2>0$ is the constant in \eqref{YQ1}.
Define
$$ I^{\le}_h := \sum_{Q\in  \mathcal W^\le_h}  h^{d-2}\L^2(F(Q)),\quad h>0;\qquad
I^{>}_j :=\sum_{Q\in \mathcal W_j} Y_Q,\quad j\in\Z.$$
For $r>0$, let $j_r=\lceil -\log_2(2C_1r)\rceil$, where $C_1$ is as in \eqref{CQW}. If $Q\in \mathcal W_j$ for some $j\in\Z$, and $F(Q)\cap (\pa D)^r\ne\varnothing$, then by \eqref{CQW}, $2^{-j-1}\le \rho_Q< C_1 r$, which implies that $j\ge j_r$. Thus,
$$ (\pa D)^r\cap D\subset\bigcup_{j=j_r}^\infty \bigcup_{Q\in\mathcal W_j} F(Q),$$
which implies that
\BGE
\begin{aligned}
 & h^{d-2}\L^2(\Gamma_D^h\cap (\pa D)^r\cap D)\le  \sum_{j=j_r}^\infty \sum_{Q\in \mathcal W_j} h^{d-2} \L^2(F(Q)\cap\Gamma_D^h)\\
 \le & \sum_{Q\in \mathcal W^\le_h} h^{d-2} \L^2(F(Q))+\sum_{j=j_r}^\infty \sum_{Q\in \mathcal W_j\cap  \mathcal W^>_h }  h^{d-2} \L^2(F(Q)\cap\Gamma_D^h) \lesssim I^\le_h+\sum_{j=j_r}^\infty I^{>}_j,
\end{aligned} \label{<I+I}
\EDE
where the last step follows from (\ref{YQ1}) since for $Q\in \mathcal W^>_h$, $h<c_2\rho_Q$.

By \eqref{CQW}, if $Q\in  \mathcal W^\le_h$, then $F(Q)\subset (\pa D)^{C_1h/c_2}$. Suppose $h<c_2 r_0/C_1$, where $r_0$ is as in \eqref{eq:boundary-q-tube}, then $C_1 h/c_2<r_0$, and \eqref{eq:boundary-q-tube} yields
\BGE  I^\le_h \le h^{d-2} \L^2( (\pa D)^{C_1h/c_2})\lesssim h^{d-2} (C_1h/c_2)^{2-q}\asymp h^{d-q}. \label{Ile}\EDE

Let $j_0=\lceil \log_2(C_1/r_0)\rceil$, and $j\ge j_0$.  By \eqref{CQW}, for every $Q\in\mathcal W_j$, $F(Q)\subset (\pa D)^{C_1 2^{-j}}$. Since $C_1 2^{-j}\le r_0$, by \eqref{eq:boundary-q-tube}, $\sum_{Q\in\mathcal W_j} \L^2(F(Q))\le \L^2 ((\pa D)^{C_1 2^{-j}})\lesssim (2^{-j})^{2-q}$.
On the other hand, for each $Q\in\mathcal W_j$, by \eqref{CQW}, $\L^2(F(Q))\gtrsim \rho_Q^2\asymp 2^{-2j}$. So we get $\#\mathcal W_j\lesssim 2^{jq}$. Since for each $Q\in\mathcal W_j$,  $\EE(Y_Q)\lesssim  \rho_Q^d\le 2^{-jd}$, we get $\EE(I^{>}_j)=\sum_{Q\in\mathcal W_j} \EE(Y_Q) \lesssim 2^{-j(d-q)}$ for $j\ge j_0$.
Since $d>q$, we get $\EE\Big(\sum_{j=j_0}^\infty  I^{>}_j\Big)=\sum_{j=j_0}^\infty \EE(I^{>}_j)<\infty$,
which further  implies a.s.\ $ \sum_{j=j_0}^\infty I_j^><\infty$.
Since $j_r\to\infty$ as $r\downarrow0$, we get a.s.\ $\lim_{r\downarrow0} \sum_{j=j_r}^\infty I_j^> =0$.
Combining this with \eqref{<I+I}, \eqref{Ile} and that $d>q$, we conclude that $\Gamma_D$ almost surely satisfies
\eqref{eq:boundary-tightness}.
\end{proof}

\begin{proposition}
Suppose that, for some $q<d$, $\lin\Cont_q(\pa D)<\infty$.
Then a.s.\ $\gamma_F$ can be reparametrized to be a curve
$\gamma_D$ defined on $[0,L]$, which is parametrized by
$d$-dimensional Minkowski content in the sense of
\eqref{eq:intervalcontent}, where
$L:=\Cont_d(\Gamma_D)$ is a random number in $(0,\infty)$.
\label{DMC}
\end{proposition}

\begin{proof}
The assumption implies that $\Cont_d(\pa D)=0$.
By Lemma~\ref{lem:target-boundary-tight},
$\Gamma_D$ almost surely satisfies
\eqref{eq:boundary-tightness}. We work on this probability-one
event. By Lemma~\ref{lem:MC-globalization} and
\eqref{eq:target-local-MC}, $ \nu :=
 F_*\bigl(
 |F'|^d\mu^\gamma_{\R_+}
 \bigr)
 =
 \Md(\cdot;\Gamma_D)$.
In particular, $L:=\nu(\Gamma_D)=\Cont_d(\Gamma_D)\in [0,\infty)$, and $\nu(\pa D)=0$.

We first record a zero-overlap property. Let
$0<r<u<\infty$. By \cite[(30)]{LawlerRezaei}, applied to
the corresponding capacity times,
$\mu^\gamma_{\R_+}(\gamma([0,u]))
 =
 \mu^\gamma_{\R_+}(\gamma([0,r]))
 +
 \mu^\gamma_{\R_+}(\gamma([r,u]))$.
It follows that $\mu^\gamma_{\R_+}( \gamma([0,r])\cap\gamma([r,u]))=0$.
By the injectivity of $F$ in $\HH$, the definition of $\nu$,
and $\nu( \pa D)=0$, we obtain
$ \nu( \gamma_F([0,r])\cap\gamma_F([r,u]))=0$.
Letting $u\uparrow\infty$ and using continuity from below gives
\begin{equation}\label{eq:adjacent-overlap-tail}
 \nu\bigl(
 \gamma_F([0,r])\cap\gamma_F([r,\infty])
 \bigr)=0,
 \qquad 0<r<\infty.
\end{equation}
The same statement holds for $r\in\{0,\infty\}$ since
$\gamma_F(\{0,\infty\})\subset\pa D$ and
$\nu(\pa D)=0$.

Define
\begin{equation}\label{eq:target-clock}
 \Theta_D (t)
 :=
 \int_0^t
 \ind_{\{\gamma(r)\in\HH\}}
 |F'(\gamma(r))|^d\,\L^1(dr),
 \qquad 0\le t\le\infty.
\end{equation}
Then $\Theta_D $ is continuous and nondecreasing, and $ \Theta_D (\infty)
 =
 \nu(\Gamma_D) = L$ by the
definition of $\nu$. By \eqref{eq:zero-boundary-time}, for almost every $r\in [0,\infty)$, $\gamma(r)\in\HH$, which implies that $\ind_{\{\gamma(r)\in\HH\}}|F'(\gamma(r))|^d>0$. Thus, $\Theta_D$ is strictly increasing.

If $0\le s\le t\le\infty$, then
$\gamma_F([0,t])=\gamma_F([0,s])\cup\gamma_F([s,t])$,
and the two sets on the right have $\nu$-null intersection.
By the definition of \(\nu\) and the preceding zero-overlap property, $\nu(\gamma_F([0,t]))=\Theta_D(t)$ and $\nu(\gamma_F([0,s]))=\Theta_D(s)$.
Therefore
\begin{equation}\label{eq:target-clock-increment}
 \nu(\gamma_F([s,t]))
 =
 \Theta_D (t)-\Theta_D (s).
\end{equation}

Define $\gamma_D := \gamma_F\circ\Theta_D ^{-1}: [0,L]\to\overline D$.
Reparametrizing \eqref{eq:adjacent-overlap-tail} gives
\begin{equation}\label{eq:adjacent-overlap-tail-D}
 \nu\bigl(
 \gamma_D([0,r])\cap\gamma_D([r,L])
 \bigr)=0,
 \qquad r\in[0,L].
\end{equation}
Let $I=[s,t]\subset[0,L]$. Since $\Gamma_D = \gamma_D([0,s]) \cup\gamma_D([s,t])\cup\gamma_D([t,L])$,
with all three sets compact,
\[
 \partial_{\Gamma_D}\gamma_D(I)
 \subset
 \bigcup_{x\in\{s,t\}}
 \left(
 \gamma_D([0,x])
 \cap
 \gamma_D([x,L])
 \right).
\]
By \eqref{eq:adjacent-overlap-tail-D}, $\nu(\partial_{\Gamma_D}\gamma_D(I))=0$.
Since $\nu=\Md(\cdot;\Gamma_D)$,
\eqref{eq:MC-continuity} and \eqref{eq:target-clock-increment} give
\[
 \Cont_d(\gamma_D(I))
 =
 \nu(\gamma_D(I))=\nu(\gamma_F([\Theta_D ^{-1}(s),\Theta_D ^{-1}(t)]))=t-s=\L^1(I).
\]
Thus $\gamma_D$ is parametrized by $d$-dimensional Minkowski
content in the sense of \eqref{eq:intervalcontent}.
\end{proof}

\section{Uniform dimension theorem for bounded domains}
\label{sec:smooth}

We now prove Theorem~\ref{thm:smooth}.  The preceding section
provides the required Minkowski-content parametrization in bounded
domains.  Indeed, if $\partial D$ is the union of finitely many
regular $C^1$ arcs, then $\pa D$ is locally connected, $\lin\Cont_1(\partial D)<\infty$,
and Proposition~\ref{DMC} applies since $1<d$.
Thus it remains to establish the H\"older and occupation estimates
needed for Corollary~\ref{cor:det}.

Throughout this section, let $D$, $a$, and $b$ be as in
Theorem~\ref{thm:smooth}. Choose a conformal map $F:\HH\to D$ whose Carath\'eodory extension,
denoted by the same symbol, satisfies $F(0)=a$ and $F(\infty)=b$. Let
$\gamma:\R_+\to\overline\HH$ be chordal SLE$_\kappa$ from $0$ to
$\infty$, parametrized by $d$-dimensional Minkowski content.
We write $\gamma_D:[0,L]\to\overline D$
for the parametrization supplied by Proposition~\ref{DMC}, which satisfies $\gamma_D\circ\Theta_D =F\circ\gamma$ for the $\Theta_D$ defined by \eqref{eq:target-clock}.


We next record the local conformal regularity needed to transfer the
estimates for $\gamma$ to $\gamma_D$ near regular $C^1$ boundary arcs.

\begin{lemma}
\label{lem:local-C1-conformal}
Let $U,U'\subset\C$ be simply connected domains and let
$G:U\to U'$ be conformal and bijective.  Let $p$ be a prime end of $U$.
Suppose that the prime-end boundaries of $U$ near $p$ and of $U'$
near $G(p)$ are represented by regular $C^1$ curves.

Then, for every $\eps,\rho\in(0,1)$, there is a sufficiently small
neighborhood $\Omega_0$ of $p$ such that $G$ is bi-$\rho$-H\"older
on $U\cap\Omega_0$, and
\begin{equation}\label{eq:C1-derivative}
 \dist(z,\pa U)^\eps
 \lesssim_{\eps,\Omega_0}
 |G'(z)|
 \lesssim_{\eps,\Omega_0}
 \dist(z,\pa U)^{-\eps},
 \qquad z\in U\cap\Omega_0.
\end{equation}
\end{lemma}

\begin{proof}
The H\"older assertion is the local form of
\cite[Corollary~1.3]{KalajC1}.

For the derivative estimate, choose $\sigma\in(0,1)$ sufficiently
close to $1$ that $ \frac1\sigma-1<\eps$ and $ 1-\sigma<\eps$.
After shrinking the neighborhoods if necessary, the distances to the
boundaries are realized on the corresponding local $C^1$ arcs.
Applying the $\sigma$-H\"older estimates to $G$ and $G^{-1}$ gives
\[
 \dist(z,\pa U)^{1/\sigma}
 \lesssim
 \dist(G(z),\pa U')
 \lesssim
 \dist(z,\pa U)^\sigma .
\]
Together with Koebe's estimate $|G'(z)|\asymp
 \frac{\dist(G(z),\pa U')}{\dist(z,\pa U)}$,
and after shrinking the neighborhood so that
$\dist(z,\pa U)<1$, this gives \eqref{eq:C1-derivative}.
\end{proof}

In applications, $G$ will be a suitable local restriction of $F$,
or of the corresponding conformal map obtained after unfolding a
positive-angle corner.

Since $\gamma$ a.s.\ visits $0$ only  at time $0$, and $\gamma(t)\ne \infty$ for $0\le t<\infty$, $\gamma_D$ visits $a$ and $b$ only at $0$ and $L$, respectively. Let $\mathcal J_D$ denote the finite set of junction points of $\partial D$.
Since $\kappa<8$, $\gamma$ a.s.\ does not visit $F^{-1}(\mathcal J_D)\sem \{0,\infty\}$; see
\cite[Theorems~6.1 and~6.4]{RohdeSchramm}.  Therefore,
\begin{equation}\label{eq:avoid-junctions}
 \PP[\gamma_D((0,L))\cap
 (\mathcal  J_D\cup\{a,b\} )
 =
 \varnothing]=1.
\end{equation}

\begin{proposition}
\label{prop:bounded-holder}
Almost surely, for every compact interval $I\subset(0,L)$ and every
$\alpha\in(0,1/d)$, the restriction $\gamma_D|_I$ is
$\alpha$-H\"older continuous.

If the endpoint condition holds at $a$, the preceding statement holds with $[0,L)$ in place of $(0,L)$. If the endpoint condition holds at $b$, the preceding statement holds with $(0,L]$ in place of $(0,L)$. If the endpoint condition holds at both $a$ and $b$,   then the preceding statement holds for $I=[0,L]$.
\end{proposition}

\begin{proof}
We work throughout on a fixed probability-one event on which
\eqref{eq:avoid-junctions}, Proposition~\ref{DMC},
Proposition~\ref{prop:chordalmass}, and
Proposition~\ref{prop:chordal-holder} hold simultaneously for all
integer time cutoffs and all rational exponents in their respective
ranges.  We fix a realization in this event.

Fix $\alpha\in(0,1/d)$.  Choose $\alpha_0\in (\alpha,1/d)\cap \Q$ and   $\rho\in(\alpha/\alpha_0,1) $. Then
$\alpha<\alpha_0\rho$. Choose  $\beta\in(1,d)\cap\Q$ and $\eps\in(0,\infty) $ such that
\begin{equation}\label{eq:choose-holder-parameters}
 \alpha<
 \frac{\alpha_0\rho}{1+\eps/(\beta-1)}.
\end{equation}

Fix a compact interval $I\subset(0,L)$.
By \eqref{eq:avoid-junctions} and the fact that the marked endpoints
$a$ and $b$ are visited only at times $0$ and $L$, respectively,
the compact set $\gamma_D(I)$ stays a positive distance from
$\mathcal J_D\cup\{b\}$.  Although $\gamma_D(I)$ is also disjoint from
$a$, we deliberately do not exclude $a$ separately here.  Indeed,
when $a\notin\mathcal J_D$, we have $F^{-1}(a)=0$, which is a finite
regular boundary point of $\HH$, and the local argument below applies
without change to compact time intervals in $[0,L)$ that contain
$0$.  In contrast, $b$ is excluded separately because
$F^{-1}(b)=\infty$, so the corresponding source preimage is not a
bounded boundary chart.

Therefore $\gamma_D(I)$ can be covered by finitely many interior
charts compactly contained in $D$ and boundary charts whose boundary
portions are regular $C^1$ arcs. Applying \eqref{eq:C1-derivative} with $\eps/d$, the boundary charts may be
chosen sufficiently small so that $F^{-1}$ is
bi-$\rho$-H\"older continuous on each of them, and
\BGE
 (\Im z)^\eps
 \lesssim_{ V}|F'(z)|^d,
 \qquad z\in F^{-1}(V\cap D). \label{derivative-bound}
\EDE
We may further assume that $\diam(F^{-1}(V))<1$ for every chart $V$.

We are going to show that $\gamma_D$ is $\alpha$-H\"older continuous on $I$, and then complete the proof of the first part. By subdividing $I$ into finitely many compact subintervals, we may
arrange that the image of each subinterval is contained in one of
these charts.  Since a finite concatenation of $\alpha$-H\"older
curves is again $\alpha$-H\"older, it suffices to prove the following uniform estimate.  For each fixed
chart $V$ and every compact subinterval $I'\subset I$ satisfying
$\gamma_D(I')\subset V$,
\begin{equation}\label{Lip-gamma-D}
 |\gamma_D(\max I')-\gamma_D(\min I')|
 \lesssim
 \L^1(I')^\alpha .
\end{equation}
Here and below until the end of the proof, all constants and other auxiliary quantities may depend on the
realization of $\gamma$, as well as on $D,a,b,F,\alpha,\alpha_0,\rho,\beta,\eps,L,I,V$, but not on $I'$, and such dependence will not be indicated.  No measurability of these quantities will be needed.

If $V$ is an interior chart, then since $\gamma$ is $\alpha_0$-H\"older continuous, $\gamma_D\circ\Theta_D=F\circ \gamma$, $\Theta_D'(t)=|F'(\gamma(t))|^d$, and $\log |F'|$ is  bounded  on compact subsets, we easily get \eqref{Lip-gamma-D}.  Suppose now $V$ is a boundary chart and let $W=F^{-1}(V)\subset \lin\HH$. 
 For $h>0$, define $ W_h:=  \{z\in W:\Im z<h\}$.
Since $W$ is bounded,   for
every $0<h\le 1$, one can find points
$z_1,\ldots,z_{N_h}\in\R$, with $N_h\lesssim h^{-1}$,
such that
\begin{equation}\label{eq:C1-strip-cover}
 W_h
 \subset
 \bigcup_{j=1}^{N_h}B(z_j,2h).
\end{equation}
Let $J =\Theta_D^{-1}(I)$ and $N_I= \lceil \max J  \rceil$. Let $I' $ be a compact subinterval of $I$ such that $\gamma_D(I')\subset V$. Let $J'=\Theta_D^{-1}(I')$. Then $J'\subset J\subset [0, N_I]$. From $F\circ\gamma=\gamma_D\circ\Theta_D$, we get $F\circ \gamma(J')= \gamma_D(I')\subset V$, and hence $\gamma(J')\subset W$.

By Proposition~\ref{prop:chordalmass}, $ \mu^\gamma_{[0,N_I]}(B(z_j,r))
 \lesssim r^{\beta}$ for any $1\le j\le N_h$ and $0<r\le 1$. Let $h\in (0,1/2)$, then from
\eqref{eq:C1-strip-cover} and $N_h\lesssim h^{-1}$,
\BGE
 \mu^\gamma_{[0,N_I]}(W_h)
 \le
 \sum_{j=1}^{N_h}
 \mu^\gamma_{[0,N_I]}(B(z_j,2h))
 \lesssim  h^{\beta-1}  .
 \label{eq:C1-strip-mass}
\EDE
Consequently, for some finite constant $C$,
\BGE
 \L^1\{x\in J':\Im \gamma(x)<h\}
 \le \mu^\gamma_{J'}(W_h)\le
 \mu^\gamma_{[0,N_I]}(W_h)
 \le
 C  h^{\beta-1}.
 \label{eq:C1-strip-time}
\EDE

Put $M:=\L^1(J)$ and  $\ell:=\L^1(J')\le M$. We may find a constant $c>0$ such that $C c^{\beta-1}<1/2$ and $c M^{1/(\beta-1)}<1/2$. Let
$ h=c \ell^{1/(\beta-1)}$, then  $C  h^{\beta-1}=Cc^{\beta-1}\ell <\ell/2$ and   $h\le c  M^{1/(\beta-1)}<1/2$, and we may apply \eqref{eq:C1-strip-time} to get $\L^1 (J'_h)\ge \ell/2$, where $J'_h:=\{x\in J':\Im \gamma(x)\ge h\}$. Since $I'=\Theta_D(J')$, \eqref{eq:target-clock} and  \eqref{derivative-bound} give
\begin{equation}\label{eq:C1-clock-lower}
\L^1(I')=\int_{J'}\ind_{\{\gamma(x)\in\HH\}} |F'(\gamma(x))|^d \L^1(dx)
 \gtrsim \int_{J'_h}  h^\eps \L^1(dx)\gtrsim  \ell^{1+\varepsilon/(\beta-1)}.
\end{equation}

On the other hand, Proposition~\ref{prop:chordal-holder} gives,  $|\gamma(\max J')-\gamma(\min J')| \lesssim \ell^{\alpha_0}$.
Since $F$ is $\rho$-H\"older  in $W$,   $\Theta_D(\max J')=\max I'$, and $\Theta_D(\min J')=\min I'$, we get
\[ |\gamma_D(\max I')-\gamma_D(\min I')|=
 |F(\gamma(\max J'))-F(\gamma(\min J'))|
 \lesssim \ell^{\alpha_0\rho}.
\]
Combining this with \eqref{eq:C1-clock-lower} and  \eqref{eq:choose-holder-parameters}, and using $\L^1(I')\le L$, we get \eqref{Lip-gamma-D}.

Suppose the endpoint condition holds at $a$.  If $a\notin\mathcal J_D$,
then the preceding argument applies without change to every compact
interval $I\subset[0,L)$.  Indeed, $F^{-1}(a)=0$ is a finite regular
boundary point of $\HH$, and the source neighborhoods appearing in
the proof above remain bounded.  Thus it remains to consider the
case $a\in\mathcal J_D$. By translation, we may assume $a=0$.

Suppose that the interior angle of $D$ at $a=0$ is
$\theta_0=q\pi$, where $q\in(0,2]$.  Since $D$ is simply connected
and $0\notin D$, there is a branch of $w^{1/q}$ in $D$, which we
denote by  $\Psi(w)=w^{1/q}$.
Choose a bounded open neighborhood $V$ of $0$ in $\C$ sufficiently
small that $\Psi$ is injective on $V\cap D$ and that
$V\cap\pa D$ consists only of the two regular $C^1$ arcs meeting
at $0$.  Let $\Phi=(\Psi|_{V\cap D})^{-1}$.
We use $\Psi$ and $\Phi$ also for their continuous extensions to the
corresponding boundary portions. The image under $\Psi$ of the two boundary arcs in $V$ forms a
regular $C^1$ arc through $0$; indeed, the $1/q$-power map unfolds
the angle $q\pi$ to $\pi$.

Choose $r>0$ sufficiently small that $U:=B(0,r)\cap\HH$
is mapped by $F$ into $V\cap D$.  Define $\ha F:=\Psi\circ F$ on $U$.
Then $\ha F$ is conformal and injective on $U$, extends
homeomorphically across the boundary portion of $U$ near $0$, and
the boundary of $\ha F(U)$ is regular $C^1$ near $\ha F(0)=0$.
Apply Lemma~\ref{lem:local-C1-conformal} to  $\ha F$ with $\eps_1/d$
at the boundary point $0$.  For every $\rho,\eps_1\in (0,1)$, there is $r'\in (0,r)$ (depending on $\rho,\eps_1$), such that for $U':=B(0,r')\cap \HH$,   $\ha F$ is
bi-$\rho$-H\"older continuous on $U'$, and
\begin{equation}\label{eq:unfolded-derivative}
  |\ha F'(z)|^d
 \lesssim(\Im z)^{-\eps_1},
 \qquad z\in U'.
\end{equation}

Choose $T_0\in(0,L)$ such that, with $u_0:=\Theta_D^{-1}(T_0)$,
we have $\gamma([0,u_0])\subset \lin {U'}$.
Define
\BGE \label{haTheta}
 \ha\Theta(u)
 :=
 \int_0^u
 \ind_{\{\gamma(s)\in\HH\}}
 |\ha F'(\gamma(s))|^d\,\L^1(ds),
 \qquad 0\le u\le u_0.
\EDE

We first check that $\ha\Theta(u_0)<\infty$.  Choose
$\beta_1\in(1,d)\cap\Q$ and
$\eps_1\in(0,\beta_1-1)$.  Put $N_0:=\lceil u_0\rceil$.
Since $U'$ is bounded, the same covering argument as in
\eqref{eq:C1-strip-cover}--\eqref{eq:C1-strip-mass} gives
\[
 \L^1\{s\in[0,u_0]:\Im\gamma(s)<h\}
 \lesssim h^{\beta_1-1},
 \qquad 0<h<1/2.
\]
Using  \eqref{eq:unfolded-derivative} and
decomposing according to dyadic values of $\Im\gamma(s)$, we obtain
\begin{align*}
 \ha\Theta(u_0)
 &\lesssim
 u_0+
 \sum_{n\ge1}
 2^{n\eps_1}
 \L^1\{s\in[0,u_0]:\Im\gamma(s)<2^{-n}\}\\
 &\lesssim
 u_0+
 \sum_{n\ge1}
 2^{-n(\beta_1-1-\eps_1)}
 <\infty.
\end{align*}
Thus, $\ha\Theta$ is absolutely continuous on $[0,u_0]$.
By \eqref{eq:zero-boundary-time}, almost every
$s\in[0,u_0]$ satisfies $\gamma(s)\in\HH$.  Since
$\ha F'$ does not vanish in $U$, it follows that $\ha\Theta$ is
strictly increasing on $[0,u_0]$.

Set $S_0:=\ha\Theta(u_0)$ and define
$\xi:[0,S_0]\to\C$ by $ \xi(\ha\Theta(u)):= \ha F(\gamma(u))$, $0\le u\le u_0$.
The regular-$C^1$ argument established above applies verbatim to
$\ha F$ and $\ha\Theta$: its only inputs are
Propositions~\ref{prop:chordalmass} and
\ref{prop:chordal-holder}, together with the local H\"older and
derivative estimates for $\ha F$ supplied by
Lemma~\ref{lem:local-C1-conformal}.
Consequently, $\xi$ is $\beta$-H\"older continuous for every
$\beta<1/d$.

Define $T=\Theta_D\circ \ha\Theta^{-1}$ on $[0,S_0]$.
Then $\gamma_D(T(s))=\Phi(\xi(s))$, $0\le s\le S_0$.
Moreover, since $F=\Phi\circ\widehat F$, by
\eqref{eq:zero-boundary-time}, the chain rule, and the definitions of
$\Theta_D$ and $\widehat\Theta$, we have
\begin{equation}\label{eq:power-clock}
 T(s_2)-T(s_1)
 =
 \int_{s_1}^{s_2}\ind_{\{\xi(s)\in \Psi(V\cap D)\}}|\Phi'(\xi(s))|^d \L^1(ds)
 =
 q^d\int_{s_1}^{s_2}
 |\xi(s)|^{d(q-1)}\L^1(ds) .
\end{equation}

Fix $\alpha<1/d$.  Choose $\eta>0$ sufficiently small that
\begin{equation}\label{eq:corner-eta-choice}
 \alpha< 1/({d+\eta/(1\wedge q)}).
\end{equation}

Fix $0\le s_1<s_2\le S_0$, and set $w_j:=\Phi(\xi(s_j))=\gamma_D(T(s_j))$, $j=1,2$, $R:=|w_1|\vee |w_2|$ and $r:=|w_2-w_1|$. If $r=0$, there is nothing to prove. Suppose $r>0$.  By the triangle inequality,  $R\ge\frac r2$. Since $V$ is bounded and $w_1,w_2\in V$, $R\lesssim 1$.
Suppose first that $|w_2|=R$.  In this case, choose $s_*$ as below
so that the terminal subinterval $[s_*,s_2]$ stays close to $w_2$.
If instead $|w_1|=R$, the same argument is applied to an initial
subinterval $[s_1,s_*]$ staying close to $w_1$.  We give the details
for the first case.

Choose $\delta\in (0,1/4)$ sufficiently small that the branch
$\psi(z)=z^{1/q}$ with $\psi(1)=1$ is defined on
$B(1,2\delta)$ and satisfies
\BGE
 |\psi(z_1)-\psi(z_2)|
 \asymp
 |z_1-z_2|,
 \qquad z_1,z_2\in B(1,2\delta). \label{varphi-lip}
\EDE

Since $w_2=\Phi(\xi(s_2))$, by continuity, there exists $s_*\in[s_1,s_2]$ such that $|\Phi(\xi(s_*))-w_2|=\delta r$
and $|\Phi(\xi(s))-w_2|\le\delta r$ for  $ s\in[s_*,s_2] $. Since $|w_2|=R$ and   $r\le2R$, we have,
\begin{equation}\label{eq:corner-annulus}
 |\xi(s)|\asymp R^{1/q},
 \qquad s\in[s_*,s_2].
\end{equation}

For $s\in[s_*,s_2]$, set $\zeta(s):= \Phi(\xi(s))/w_2$.
Then $\zeta(s_2)=1$, and, since $r\le2R$,
\[
 |\zeta(s)-1|
 =
 {|\Phi(\xi(s))-w_2|}/{R}
 \le  {\delta r}/{R}
 \le 2\delta ,
 \qquad s\in[s_*,s_2].
\]

Since  $\Psi(z)$ is a branch of $z^{1/q}$, the same holds for $\Psi_{w_2}(z):=\Psi(w_2z)/\Psi(w_2)$. Since $\psi$ is also a branch of $z^{1/q}$, $\psi(1)=1=\Psi_{w_2}(1)$, and $\zeta(s_2)=1$, we get by the continuity of $\zeta$,
$$\psi(\zeta(s))=\Psi_{w_2}(\zeta(s))=\Psi(w_2\zeta(s))/\Psi(w_2)=\xi(s)/\xi(s_2),\quad s\in [s_*,s_2].$$
Since $|\xi(s_2)|=|\Psi(w_2)|=R^{1/q}$, \eqref{varphi-lip} yields
\begin{align}
 |\xi(s_2)-\xi(s_*)|
 &=
 R^{1/q}|\psi(\zeta(s_2))-\psi(\zeta(s_*))|
 \asymp
 R^{1/q}|\zeta(s_2)-\zeta(s_*)| \notag\\
 &=
 R^{1/q-1} {|w_2-\Phi(\xi(s_*))|}
 =
 \delta rR^{1/q-1}
 \asymp rR^{1/q-1}.
 \label{eq:corner-source-displacement}
\end{align}

Since $1/(d+\eta)<1/d$, $\xi$ is $1/(d+\eta)$-H\"older continuous,
\eqref{eq:corner-source-displacement} gives
\begin{equation}\label{eq:corner-source-time}
 s_2-s_*
 \gtrsim
 \bigl(r R^{1/q-1}\bigr)^{d+\eta}.
\end{equation}

Using \eqref{eq:power-clock}, \eqref{eq:corner-annulus}, and
\eqref{eq:corner-source-time}, we obtain
\begin{align}
 T(s_2)-T(s_1)
 &\ge
 q^d\int_{s_*}^{s_2}
 |\xi(s)|^{d(q-1)}\,\L^1(ds) \gtrsim
 R^{d(1-1/q)}(s_2-s_*) \notag\\
 &\gtrsim
 R^{d(1-1/q)}
 \bigl(r R^{1/q-1}\bigr)^{d+\eta}=
 r^{d+\eta}R^{(1/q-1)\eta}.
 \label{eq:corner-duration}
\end{align}

Since   $r/2\le R\lesssim 1$, $R^{(1/q-1)\eta}\gtrsim  r^{0\vee (1/q-1)\eta}$. From \eqref{eq:corner-duration} we then get $T(s_2)-T(s_1)\gtrsim r^{d+(1\vee (1/q))\eta}$,
or equivalently
\begin{equation}\label{eq:corner-holder}
|w_1-w_2|= r\lesssim
 |T(s_2)-T(s_1)|^{1/(d+\eta/(1\wedge q)}.
\end{equation}

Since $w_j=\gamma_D(T(s_j))$, by  \eqref{eq:corner-eta-choice} and
\eqref{eq:corner-holder}, we get, for any $0\le s_1\le s_2\le S_0$,
\[
 |\gamma_D(T(s_2))-\gamma_D(T(s_1))|
 \lesssim
 |T(s_2)-T(s_1)|^\alpha.
\]
Since $T$ maps $[0,S_0]$ onto $[0,T_0]$, this proves that
$\gamma_D$ is $\alpha$-H\"older continuous on $[0,T_0]$.

Finally, for every $T<L$, the first part of the proof gives the
$\alpha$-H\"older estimate on every compact subinterval of $(0,T]$,
while the preceding argument gives it on a neighborhood of $0$.
A finite-cover argument therefore shows that $\gamma_D$ is
$\alpha$-H\"older continuous on $[0,T]$.

Suppose the endpoint condition holds at $b$.  Define the reversed
curve $\widetilde\gamma_D(t):=\gamma_D(L-t)$, $0\le t\le L$.
Since $\gamma_D$ is parametrized by $d$-dimensional Minkowski
content, so is $\widetilde\gamma_D$. By reversibility of chordal SLE$_\kappa$
\cite{ZhanReversibility,MSIG3}, $\widetilde\gamma_D$ has the law of
chordal SLE$_\kappa$ in $D$ from $b$ to $a$, parametrized by
$d$-dimensional Minkowski content.  Applying the initial-endpoint
assertion proved above to $\widetilde\gamma_D$, we conclude that  almost surely, for every compact
interval $I\subset(0,L]$, $\gamma_D|_I$ is $\alpha$-H\"older continuous.  This proves the third assertion.

Finally, suppose that the endpoint condition holds at both
$a$ and $b$.  The preceding two assertions imply that, for any $\alpha\in (0,1/d)$,
$\gamma_D$ is $\alpha$-H\"older continuous on both $[0,L/2]$ and $[L/2,L]$.
It follows that $\gamma_D$ is $\alpha$-H\"older continuous on
$[0,L]$.
\end{proof}

\begin{proposition}
\label{prop:bounded-occupation}
There exists an event of probability one on which, for every compact
interval $I\subset(0,L)$ and every $0<\beta<d$, there is a finite random
constant $C_{I,\beta}(\gamma)$ such that
\begin{equation}\label{eq:bounded-local-mass}
 \mu^{\gamma_D}_I(B(z,r))
 \le C_{I,\beta}(\gamma) r^\beta,
 \qquad z\in\C,\quad 0<r\le1.
\end{equation}

If the endpoint condition holds at $a$, then  \eqref{eq:bounded-local-mass} holds for any compact
interval $I\subset[0,L)$. If the endpoint condition holds at $b$, then  \eqref{eq:bounded-local-mass} holds for any compact
interval $I\subset(0,L]$.  If the endpoint condition holds at both $a$ and
$b$, then \eqref{eq:bounded-local-mass} holds for $I=[0,L]$.
\end{proposition}

\begin{proof}
We work throughout on the same probability-one event as in the proof
of Proposition~\ref{prop:bounded-holder}, and fix a realization in
this event.

Fix a compact interval $I\subset(0,L)$ and an exponent
$\beta\in(0,d)$.  Choose $\rho\in(0,1)$, $\beta_0\in(1\vee \beta,d)\cap \Q$, and $\eps\in (0,\beta_0-1)$
so that
\begin{equation}\label{eq:occupation-parameters}
 \beta<\rho(\beta_0-\eps).
\end{equation}

As in the proof of Proposition~\ref{prop:bounded-holder},
$\gamma_D(I)$ can be covered by finitely many interior charts and
regular $C^1$ boundary charts.  We choose, for each such chart $V$, a
smaller chart $V_0$ whose closure is contained in  $V$ on
which the corresponding conformal estimates hold.  It suffices to
estimate the contribution from one such $V_0$.

Here and below, all constants and auxiliary quantities may depend on
the realization of $\gamma$ and on
$D,a,b,F,I,\beta,\rho,\beta_0,\eps$ and the charts fixed above, but not
on the center or radius of the ball under consideration.  Such dependence will not be indicated. For
convenience of exposition, we may shrink  these
charts further below; any such adjustment depends only on the data
already fixed at this stage and could equivalently have been made at
the outset.

Let $J=\Theta_D^{-1}(I)$ and
$N_I=\lceil\max J\rceil$.  By \eqref{eq:target-clock}, for every
Borel set $E\subset D$,
\begin{equation}\label{eq:weighted-occupation-bounded}
 \mu^{\gamma_D}_I(E)
 =
 \int_J
 \ind_{\{F(\gamma(t))\in E\}}
 |F'(\gamma(t))|^d\,\L^1(dt).
\end{equation}

Suppose first that $V$ is an interior chart.  Since
$ V_0\Subset V\Subset D$, both $F^{-1}$ and $F'$ are
bounded and Lipschitz on the relevant compact sets.  Hence there
exists $r_0>0$ such that, whenever $0<r\le r_0$ and $B(w,r)$ meets
$\gamma_D(I)\cap V_0$, its $F$-preimage is contained in a   ball
of radius $O(r)\le1$.  Proposition~\ref{prop:chordalmass} and \eqref{eq:weighted-occupation-bounded} therefore
give $\mu^{\gamma_D}_I(B(w,r)\cap V_0)\lesssim r^{\beta_0}
 \lesssim r^\beta$ for $0<r\le r_0$.  If $r_0\le r\le1$, the same bound follows from $\mu^{\gamma_D}_I(B(w,r)\cap V_0)\le \L^1(I)\le L\lesssim 1$.

Suppose now that $V$ is a boundary chart, and put $W:=F^{-1}(V\cap D)$.
By Lemma~\ref{lem:local-C1-conformal},  $W$ is bounded, $F^{-1}$ is $\rho$-H\"older continuous on $V$, and
\begin{equation}\label{eq:occupation-derivative-bound}
 |F'(z)|^d\lesssim(\Im z)^{-\eps},
 \qquad z\in W.
\end{equation}

Since $\overline V_0\subset V$, we may choose $r_V>0$ so small that $(\lin V_0)^{2r_V}\cap\lin D\subset V$
and also $2^\rho C_\rho r_V^\rho\le1$, where $C_\rho$ is the $\rho$-H\"older constant of
$F^{-1}$ on $V$. Let $0<r\le r_V$ and $\til r=2^\rho C_\rho r^\rho\le 1$. If $ \mu^{\gamma_D}_I(B(w,r)\cap V_0)>0$, choose $x\in B(w,r)\cap V_0\cap\gamma_D(I)$.  Since
$B(w,r)\cap \lin D\subset V_0^{2r}\cap \lin D \subset V$,
the $\rho$-H\"older continuity of $F^{-1}$ on $V$ gives
\BGE F^{-1}(B(w,r)\cap \lin D)
 \subset B(F^{-1}(x),\til r)\cap W.\label{FBwr}\EDE

For $0<h\le \til r$, the set $B(F^{-1}(x),\til r)\cap W\cap\{\Im z<h\}$
can be covered by $O(\til r /h)$ balls of radius $h$.  Hence
Proposition~\ref{prop:chordalmass} gives
\begin{equation}\label{eq:occupation-strip}
 \mu^\gamma_{[0,N_I]}
 \bigl(B(F^{-1}(x),\til r)\cap W\cap\{\Im z<h\}\bigr)
 \lesssim \til r h^{\beta_0-1}.
\end{equation}

We now decompose $B(F^{-1}(x),\til r)\cap W$ according to the distance from
$\R$.  On the part where $\Im z\ge \til r$, using
\eqref{eq:occupation-derivative-bound} and
Proposition~\ref{prop:chordalmass},
\[
 \int_{B(F^{-1}(x),\til r)\cap W\cap\{\Im z\ge \til r\}}
 |F'(z)|^d\,\mu^\gamma_{[0,N_I]}(dz)
 \lesssim
 \til r^{-\eps}\til r^{\beta_0}
 =
 \til r ^{\beta_0-\eps}.
\]
On the dyadic strip $S_n:=\{z\in\C:2^{-n}\til r\le\Im z<2^{-n+1}\til r\}$,
\eqref{eq:occupation-derivative-bound} and
\eqref{eq:occupation-strip} give
\[
 \int_{S_n\cap B(F^{-1}(x),\til r)\cap W}
 |F'(z)|^d\,\mu^\gamma_{[0,N_I]}(dz)
 \lesssim
 \til r(2^{-n}\til r)^{\beta_0-1-\eps},\quad n\in\N.
\]
Since $\eps<\beta_0-1$, summing the above two displayed formulas and using \eqref{FBwr} and \eqref{eq:occupation-parameters}  yields
\[
 \mu^{\gamma_D}_I(B(w,r))=\mu^{\gamma_D}_I(B(w,r)\cap D)
 \lesssim
\til r^{\beta_0-\eps}\asymp r^{\rho(\beta_0-\eps)}
 \lesssim
 r^\beta.
\]
Here we used $\mu^{\gamma_D}(\pa D)=0$. For $r_V<r\le1$, the same estimate follows after increasing the
constant, since $\mu^{\gamma_D}_I(\C)=\L^1(I)\le L\lesssim 1$.  Summing over
the finite chart cover proves \eqref{eq:bounded-local-mass}.

We next assume that the endpoint condition holds at $a$.  If $a$ is
not a junction point, the preceding argument applies without change
to every compact interval $I\subset[0,L)$, since we may still cover $\gamma_D(I)$ by finitely many interior charts and regular $C^1$ boundary charts, and the charts can be chosen such that their $F$-preimages are still bounded.

Suppose now that $a$ is a junction point with positive interior
angle $q\pi$, where $q\in(0,2]$.  We briefly recall the unfolding
construction from the proof of Proposition~\ref{prop:bounded-holder}.
After translating, we may assume that $a=0$.  Let $ \Psi(w)$ be a branch of $w^{1/q}$ in $D$.  For a sufficiently small bounded
neighborhood $V$ of $0$ in $\lin D$, the map $\Psi$ is injective on
$V\cap D$, and the two boundary arcs of $D$ meeting at $0$ are
unfolded by $\Psi$ into a regular $C^1$ boundary arc.  Write
 $\Phi=(\Psi|_{V\cap D})^{-1}$.

Choose $T_0\in(0,L)$ sufficiently small that $\gamma_D([0,T_0])\subset V$,
and put $u_0:=\Theta_D^{-1}(T_0)$.
On $F^{-1}(V)$ define the unfolded conformal map $\ha F:=\Psi\circ F$
and its natural clock $\ha\Theta$ on $[0,u_0]$ by \eqref{haTheta}, which satisfies
$\ha\Theta(u_0)<\infty$.  Set
 $S_0:=\ha\Theta(u_0)$,
and define the unfolded curve $\xi:[0,S_0]\to\C$ by $\xi=\ha F\circ \gamma\circ \ha\Theta^{-1}$. Then
 $\xi$ is   parametrized by  $d$-dimensional Minkowski
content.

Define $T:=\Theta_D\circ\ha\Theta^{-1}:[0,S_0]\to[0,T_0]$. We have
$ \gamma_D(T(s)) = \Phi(\xi(s))$, $0\le s\le S_0$.
Moreover, we have \eqref{eq:power-clock}. Thus, for any Borel set $S\subset \C$,
\BGE
 \mu^{\gamma_D}_{[0,T_0]}(S)=\mu^{\gamma_D}_{[0,T_0]}(S\cap V\cap D)=q^d \int_{\Psi(S\cap V\cap D)} |w|^{d(q-1)}\,\mu^\xi_{[0,S_0]}(dw).\label{eq:occupation-power-clock}
 \EDE

Let $\eps_1=\min\{dq,(1\wedge q)(d-\beta)\}/2$. Then
\begin{equation}\label{eq:corner-occupation-eps}
 dq-\eps_1>0,
 \qquad
 \beta< d- {\eps_1}/({1\wedge q}).
\end{equation}
The target boundary for $\ha F$ is regular $C^1$ near $0$.
Therefore the regular-$C^1$ occupation estimate proved above,
applied to $\ha F$ and the clock $\ha\Theta$, yields, for $0<r\le 1$,
\begin{equation}\label{eq:unfolded-occupation}
 \mu^\xi_{[0,S_0]}(B(z,r))
 \lesssim
 r^{d-\eps_1},
 \qquad z\in\C.
\end{equation}
After increasing the implicit constant, \eqref{eq:unfolded-occupation} holds for
all $r>0$, since $\mu^\xi_{[0,S_0]}(\C)=S_0<\infty$.

Since we have proved \eqref{eq:bounded-local-mass} for compact subintervals of $(0,L)$, it now suffices to show that \eqref{eq:bounded-local-mass} holds for $[0,T_0]$.
Let $C_q=1/\sin(q\pi/12)$. Then $C_q\ge 2$ since $0<q\le 2$. Suppose first that $|z|\le C_q r$. Let $R_0:=((C_q+1)r)^{1/q}$. Since $\Psi(w)$ is a branch of $w^{1/q}$,  $ \Psi(B(z,r)\cap V\cap D) \subset B(0,R_0)$.  Decomposing $B(0,R_0)$ into dyadic annuli at scales
$2^{-n}R_0$, $n\in\N\cup\{0\}$, and using
\eqref{eq:occupation-power-clock} and
\eqref{eq:unfolded-occupation}, we obtain
$$ \mu^{\gamma_D}_{[0,T_0]}(B(z,r))
 \lesssim
 \sum_{n\ge0}
 (2^{-n}R_0)^{d(q-1)} (2^{-n}R_0)^{d-\eps_1}
 \lesssim
 R_0^{dq-\eps_1}
 \asymp
 r^{d-\eps_1/q}\le r^\beta,$$
where we used \eqref{eq:corner-occupation-eps} in the last and the third-to-last inequalities.

Suppose next that $|z|>C_qr$.  We claim that
$\Psi(B(z,r)\cap V\cap D)$ can be covered by at most two balls of
radius $ K_q r|z|^{1/q-1}$, where $K_q:= 2^{|1/q-1|}/q$.

Let $\vartheta(w)$ be the continuous branch of $\arg w$ in $D$ so that  $ \Psi(w)=|w|^{1/q}e^{i\vartheta(w)/q}$.
By choosing $V$ sufficiently small, we may assume that  $ \sup_{ V\cap D}\vartheta - \inf_{ V\cap D}\vartheta <3\pi$.

Since $0\notin B(z,r)$, there is a continuous branch $\vartheta_B(w)$ of $\arg w$ on $B(z,r)$. Define
 $\psi_B(w):=|w|^{1/q}e^{i\vartheta_B(w)/q}$. Then $\psi_B(w)$ is a branch of $w^{1/q}$ in $B(z,r)$.
We observe  that
\[
 \sup_{B(z,r)}\vartheta_B-  \inf_{B(z,r)}\vartheta_B= 2\arcsin(r/|z|)
 < 2\arcsin(C_q^{-1})
=\frac{q\pi}{6}.
\]
For $w\in B(z,r)\cap V\cap D$, $(\vartheta(w)-\vartheta_B(w))/(2\pi)\in\Z$.
At most two integers can occur.  Indeed, if three distinct integers
occurred, the largest and smallest would differ by at least $2$, and we would have
\[
 \sup_{V\cap D}\vartheta- \inf_{V\cap D}\vartheta
 \ge
 4\pi-(\sup_{B(z,r)}\vartheta_B- \inf_{B(z,r)}\vartheta_B)
 >
 4\pi-\frac{q\pi}{6}
 \ge\frac{11\pi}{3},
\]
contradicting that  $ \sup_{ V\cap D}\vartheta - \inf_{ V\cap D}\vartheta <3\pi$.

For each of the at most two integers $k$ that occurs, on the
corresponding subset of $B(z,r)\cap V\cap D$ we have
$\Psi(w)=e^{2\pi i k/q}\psi_B(w)$.
Moreover, since $C_q\ge 2$, $ \frac12|z|\le |w|\le 2|z|$ for $w\in B(z,r)$,
and therefore $ |\psi_B'(w)|
 = |w|^{1/q-1}/q
 \le K_q |z|^{1/q-1}$.
As $B(z,r)$ is convex,
\[
 |\psi_B(w)-\psi_B(z)|
 \le K_q r|z|^{1/q-1},
 \qquad w\in B(z,r).
\]
Thus each of the at most two subsets is mapped by $\Psi$ into one
ball of radius $K_qr|z|^{1/q-1}$, proving the claim.

Let $R=\sup_{\zeta\in V}|\zeta|<\infty$. If $B(z,r)\cap V\ne\varnothing$, then $|z|\le R+r\le R+1$. If $1/q-1\ge 0$, then $r|z|^{1/q-1}\le (R+1)^{1/q-1}r$. If $1/q-1\le 0$, then from $|z|\ge C_q r$ we get $r|z|^{1/q-1}\le r(C_q r)^{1/q-1}\lesssim r^{1/q}$. Both upper bounds tend to $0$ as $r\downarrow 0$, implying that we can find a constant $r_1\in (0,1]$ such that $K_q r|z|^{1/q-1}\le 1$ whenever $0<r\le r_1$ and $B(z,r)\cap V\ne\varnothing$.

We now bound $ \mu^{\gamma_D}_{[0,T_0]}(B(z,r))$. If $B(z,r)\cap V=\varnothing$, then since $\gamma_D([0,T_0])\subset V$, we have $ \mu^{\gamma_D}_{[0,T_0]}(B(z,r))=0$. Suppose now $B(z,r)\cap V\ne \varnothing$ and $0<r\le r_1$. By \eqref{eq:occupation-power-clock},
\eqref{eq:unfolded-occupation}, and the facts that  $|w|\asymp |z|^{1/q}$ for
$w\in \Psi(B(z,r)\cap V\cap D)$, and $\Psi(B(z,r)\cap V\cap D)$ is covered by at most two balls of radius $K_qr|z|^{1/q-1}\le 1$, we get
\begin{align}
 \mu^{\gamma_D}_{[0,T_0]}(B(z,r))
 &\lesssim
 |z|^{d(1-1/q)}
 \mu^\xi_{[0,S_0]}
 \bigl(\Psi(B(z,r)\cap V\cap D)\bigr)
 \notag\\
 &\lesssim
 2 |z|^{d(1-1/q)}
 (K_q r|z|^{1/q-1})^{d-\eps_1}\asymp
 r^{d-\eps_1}|z|^{(1-1/q)\eps_1}. \notag
\end{align}

If $1-1/q\ge 0$, then from $|z|\le R+1$, we get
$ \mu^{\gamma_D}_{[0,T_0]}(B(z,r))
 \lesssim r^{d-\eps_1}\le r^\beta$.
If $1-1/q\le 0$, then from $|z|\ge C_q r$, we get
$ \mu^{\gamma_D}_{[0,T_0]}(B(z,r))
 \lesssim
r^{d-\eps_1}r^{(1-1/q)\eps_1}=r^{d-\eps_1/q}\le r^\beta$.
Thus, in both cases, $ \mu^{\gamma_D}_{[0,T_0]}(B(z,r))\lesssim r^\beta$ if $0<r\le r_1$. The inequality also holds if $r\ge r_1$ since $ \mu^{\gamma_D}_{[0,T_0]}(B(z,r))\le L\lesssim 1$.
Thus \eqref{eq:bounded-local-mass} holds with $I=[0,T_0]$.

Together with the already proved estimate on compact subintervals of
$(0,L)$, this proves \eqref{eq:bounded-local-mass} for every compact
interval $I\subset[0,L)$ whenever the endpoint condition holds at
$a$.

The terminal endpoint is treated exactly as in the proof of
Proposition~\ref{prop:bounded-holder}.  If the endpoint condition
holds at $b$, reverse the curve by setting
$\widetilde\gamma_D(t)=\gamma_D(L-t)$.  By reversibility of chordal
SLE$_\kappa$ \cite{ZhanReversibility,MSIG3} and the orientation
invariance of Minkowski content, the initial-endpoint assertion
proved above gives \eqref{eq:bounded-local-mass} for every compact
interval $I\subset(0,L]$.

Finally, if the endpoint condition holds at both $a$ and $b$, the
preceding two assertions applied to $[0,L/2]$ and $[L/2,L]$,
respectively, and the additivity of the occupation measure give
\eqref{eq:bounded-local-mass} for $I=[0,L]$.
\end{proof}

\begin{proof}[Proof of Theorem~\ref{thm:smooth}]
As observed at the beginning of this section,
$\lin\Cont_1(\partial D)<\infty$ and $1<d$.
Proposition~\ref{DMC} therefore gives a parametrization
$ \gamma_D:[0,L]\to\overline D $
by $d$-dimensional Minkowski content in the sense of \eqref{eq:intervalcontent}, with
$0<L<\infty$ almost surely.
Intersect the probability-one events in
Propositions~\ref{prop:bounded-holder} and
\ref{prop:bounded-occupation}.  Fix a path in this common event.

For $n\ge3$, put $I_n= [\frac{L}{n}, (1-\frac1n )L ]$.
On $I_n$, Propositions~\ref{prop:bounded-holder} and
\ref{prop:bounded-occupation} give every H\"older exponent below
$1/d$ and every occupation exponent below $d$.  Hence
Corollary~\ref{cor:det} applies to $\gamma_D|_{I_n}$ and gives
\eqref{eq:main} with $\gamma_D$ in place of $\eta$, simultaneously for every $A\subset I_n$.

For arbitrary $A\subset[0,L]$, $A\cap(0,L)
 = \bigcup_{n\ge3}(A\cap I_n)$.
Hausdorff and packing dimensions are countably stable, and adding the
two endpoint times does not change either dimension.  It follows that,
simultaneously for every $A\subset[0,L]$, \eqref{eq:main} holds for $\dimH$ and $\dimP$ with $\gamma_D$ in place of $\eta$.

Suppose now that the endpoint condition holds at both $a$ and $b$.
Propositions~\ref{prop:bounded-holder} and
\ref{prop:bounded-occupation} then give, on the full interval $[0,L]$,
every H\"older exponent below $1/d$ and every occupation exponent
below $d$.  Applying Corollary~\ref{cor:det} directly to
$f=\gamma_D$ gives \eqref{eq:main} with $\gamma_D$ in place of $\eta$ for $\dimBl$ and $\dimBu$ and every $A\subset [0,L]$.
This completes the proof.
\end{proof}

\begin{remark} \label{rem:simple-bounded-domain}
Suppose that $0<\kappa\le4$.  In this case the global boundary
regularity assumptions above can be localized entirely to the two
marked endpoints.  To see this, let $D$ be a bounded simply connected
domain and let $F:\HH\to D$ be conformal.  Suppose that $F$ extends
continuously at $0$ and $\infty$, with $F(0)=a$ and $F(\infty)=b$
and that, for some $q<d$ and $r_a,r_b>0$,
\[
 \lin\Cont_q\bigl(\pa D\cap\lin B(a,r_a)\bigr)
 +
 \lin\Cont_q\bigl(\pa D\cap\lin B(b,r_b)\bigr)
 <\infty.
\]
Let $\gamma$ be chordal SLE$_\kappa$ in $\HH$ from $0$ to $\infty$,
parametrized by $d$-dimensional Minkowski content, and set $
 \gamma_F(t)=F(\gamma(t))$, $0\le t\le \infty$ with $\gamma_F(0)=a$ and $\gamma_F(\infty)=b$.   Write
$ \Gamma_D:=\gamma_F([0,\infty])$.
Almost surely, $ \gamma((0,\infty))\subset\HH$,
and hence
 $\gamma_F((0,\infty))\subset D$.
By conformal covariance of Minkowski content measure,
$ \nu:= F_*\bigl(|F'|^d\mu^\gamma_{\R_+}\bigr)=\Md(\cdot;\Gamma_D\cap D)$.

Moreover, by continuity at the two marked endpoints, there exist
$0<u<v<\infty$ such that $\gamma_F([0,u])\subset B(a,r_a)$ and
 $\gamma_F([v,\infty])\subset B(b,r_b)$,
while $\gamma_F([u,v])\Subset D$.  Hence the boundary-tightness
argument in Proposition~\ref{DMC} need only be applied to the two
endpoint tails, using the corresponding local bounds on $\pa D$.
The globalization argument localizes in the same way.  Since
$\Gamma_D\cap\pa D=\{a,b\}$ and the local $q$-dimensional Minkowski
bounds imply zero $d$-dimensional Minkowski content near these points,
we obtain almost surely $\nu=\Md(\cdot;\Gamma_D)$ and
 $\Cont_d(\Gamma_D)=\nu(\Gamma_D)\in(0,\infty)$.

Define $\Theta_D$ using \eqref{eq:target-clock}.
Then $\Theta_D$ is an increasing homeomorphism from $[0,\infty]$ onto $[0,L]$ for $L=\nu(\Gamma_D)\in(0,\infty)$, and the curve  $\gamma_D:= \gamma_F\circ\Theta_D^{-1}$ defined on $[0,L]$
is parametrized by $d$-dimensional Minkowski content, including the
two endpoints.

Once this finite Minkowski-content parametrization is available,
the conclusions of Propositions~\ref{prop:bounded-holder} and
\ref{prop:bounded-occupation} on compact intervals
$I\subset(0,L)$ require no boundary regularity at all: since
$ \gamma_D(I)\Subset D$, only interior charts are needed.  To include $0$ or $L$ in the
time interval, it suffices to impose the local piecewise regular
$C^1$ assumption and the endpoint condition near the corresponding
endpoint.

For the Hausdorff and packing dimension conclusions, the endpoint
assumptions can be removed altogether.  Indeed, redefine $\Theta_D$ by
$ \Theta_D(t):=\int_1^t |F'(\gamma(s))|^d\L^1(ds)$, $0<t<\infty$.
Then $\Theta_D$ is an increasing homeomorphism from $(0,\infty)$
onto an open interval $J\subset\R$ containing $0$.  Defining $\gamma_D$ on $J$ in the same way, we obtain a parametrization by $d$-dimensional Minkowski content.
Moreover, on a single probability-one event, for every $A\subset J$,
\eqref{eq:main} holds for $\gamma_D|_J$ and
$\mathfrak D\in\{\dimH,\dimP\}$.  Indeed, one first proves the claim
for $A\subset I\Subset J$, where only interior charts are needed, and
then uses a countable exhaustion of $J$ together with the countable
stability of Hausdorff and packing dimensions.
\end{remark}

\section{Exact local dimension of natural occupation}\label{sec:localdim}

The two pathwise inputs also give a pointwise statement which is stronger than the
usual almost-everywhere local dimension of a measure.

\begin{proposition}
\label{prop:localdim}
Let $(\gamma,I)$ be one of the following:
\begin{enumerate}
\item $\gamma=\eta$ is the sssi SLE$_\kappa$ of
Theorem~\ref{thm:main}, and $I\subset\R$ is a nondegenerate compact
interval;

\item $\gamma$ is the chordal SLE$_\kappa$ of
Theorem~\ref{thm:chordal}, and $I\subset[0,\infty)$ is a
nondegenerate compact interval;

\item $\gamma=\gamma_D$ is the bounded-domain curve of
Theorem~\ref{thm:smooth}, and $I\subset(0,L)$ is a nondegenerate
compact interval.
\end{enumerate}
Then, on a probability-one event,
\[ \lim_{r\downarrow0}
 \frac{\log\mu_I^\gamma(B(z,r))}{\log r}
 =d,\qquad \text{uniformly in }z\in\gamma(I).
 \]
Equivalently, for every $\varepsilon\in (0,d)$ there are finite positive
random constants $c,C,r_0$ such that
\begin{equation}\label{eq:twosidedmass}
 c r^{d+\varepsilon}
 \le \mu_I^\gamma(B(z,r))
 \le C r^{d-\varepsilon},
 \qquad
 z\in\gamma(I),\quad 0<r<r_0.
\end{equation}
In case {\rm(3)}, $I$ may include $0$ if the endpoint condition holds
at $a$, may include $L$ if it holds at $b$, and may equal $[0,L]$ if
it holds at both endpoints.
\end{proposition}
\begin{proof}
The upper bound of \eqref{eq:twosidedmass} is given by the occupation power bound in the corresponding proposition among Propositions~\ref{prop:sssi-occupation},
\ref{prop:chordalmass}, and
\ref{prop:bounded-occupation}.

For the lower bound, choose $\alpha=\frac1{d+\varepsilon}<\frac1d$.
The H\"older continuity in  the corresponding proposition among Propositions~\ref{prop:sssi-holder}, \ref{prop:chordal-holder}, and \ref{prop:bounded-holder} gives an almost surely finite positive $\alpha$-H\"older constant $H$ of $\gamma$ on $I$. Let $r_0=H(\L^1(I)/2)^\alpha\wedge 1$. If $0<r<r_0$, then $(r/H)^{1/\alpha}<\L^1(I)/2$.  If $z=\gamma(t)$ with $t\in I$, then there exists a subinterval
$I'\subset I$ such that  $\L^1(I')=(r/(2H))^{1/\alpha}$
and $t$ is one of the endpoints of $I'$.  By the H\"older continuity of $\gamma$ on $I$, $I'$ is mapped into $B(z,r)$.  Hence
$ \mu_I^\gamma(B(z,r))\ge \L^1(I') =(2H)^{-(d+\eps)} r^{d+\eps}$.
This proves \eqref{eq:twosidedmass} in all three cases.
\end{proof}

\begin{remark}\label{rem:log-lower}
In the sssi case, the lower bound in
Proposition~\ref{prop:localdim} can be sharpened.
By the sharp modulus of continuity
\cite[Theorem~1.2]{HoldenYuan}, on a probability-one event, for every
nondegenerate compact interval $I\subset\R$, there are finite positive
random constants $c$ and $r_0$ such that
\[
 \mu_I^\eta(B(z,r))
 \ge c r^d(\log(1/r))^{1-d},
 \qquad z\in\eta(I),\quad 0<r<r_0.
\]
\end{remark}

\section*{Use of artificial intelligence}
The author used GPT~6 to explore proof ideas, locate relevant
references, and assist with the exposition.
The author simplified the proofs, reorganized the arguments,
revised the exposition, carefully checked the mathematical
details, and reviewed the final manuscript.
The author takes full responsibility for its contents.

\end{document}